\documentclass[aap,seceqn,citesort,MSNbibl,dvips]{arximspdf}
\usepackage{graphicx}
\doi{10.1214/10-AAP756}
\volume{21}
\issue{6}
\pubyear{2011}
\firstpage{2226}
\lastpage{2262}

\makeatletter

\newtheorem{theorem}{Theorem}[section]
\newtheorem{lemma}[theorem]{Lemma}

\newproclaim{remark}{Remark}
\newproclaim{algorithm}{Algorithm}
\newproclaim{example}[theorem]{Example}

\newcommand{\R}{\mathbb R}
\newcommand{\Z}{\mathbb Z}
\newcommand{\LL}{\mathbb{L}^V}
\newcommand{\E}{\mathbb{E}}

\newcommand{\eqdef}{\stackrel{\mathrm{def}}{=}}

\newcommand{\MT}{\mathcal{Z}}
\newcommand{\Err}{\mathcal{E}}

\newcommand{\EG}{\mathcal{A}}
\newcommand{\TG}{\mathcal{B}}
\newcommand{\EGV}{\mathcal{A}^V}
\newcommand{\TGV}{\mathcal{B}^V}

\makeatother

\begin{document}
\begin{frontmatter}

\title{Error analysis of tau-leap simulation methods}
\runtitle{Error analysis of tau-leap methods}

\begin{aug}
\author[A]{\fnms{David F.} \snm{Anderson}\corref{}\thanksref{t1}\ead[label=e1]{anderson@math.wisc.edu}\ead[url,label=u1]{http://www.math.wisc.edu/\textasciitilde anderson/}},
\author[B]{\fnms{Arnab} \snm{Ganguly}\thanksref{t2}\ead[label=e2]{gangulya@control.ee.ethz.ch}\ead[url,label=u2]{http://control.ee.ethz.ch/\textasciitilde gangulya/}} and
\author[A]{\fnms{Thomas G.} \snm{Kurtz}\thanksref{t2}\ead[label=e3]{kurtz@math.wisc.edu}\ead[url,label=u3]{http://www.math.wisc.edu/\textasciitilde kurtz/}}
\runauthor{D. F. Anderson, A. Ganguly and T. G. Kurtz}
\affiliation{University of Wisconsin, Madison}
\address[A]{D. F. Anderson\\
T. G. Kurtz\\
Department of Mathematics\\
University of Wisconsin\\
480 Lincoln Drive\\
Madison, Wisconsin 53706\\
USA\\
\printead{e1}\\
\phantom{E-mail: }\printead*{e3}\\
\printead{u1}\\
\phantom{URL: }\printead*{u3}}
\address[B]{A. Ganguly\\
Automatic Control Laboratory\\
Swiss Federal Institute of Technology\\
Physikstrasse 3, ETL I13\\
8092 Z\"{u}rich\\
Switzerland\\
\printead{e2}\\
\printead{u2}}
\end{aug}

\thankstext{t1}{Supported by NSF Grants DMS-05-53687 and DMS-10-09275.}
\thankstext{t2}{Supported by NSF Grant DMS-05-53687.}

\received{\smonth{9} \syear{2009}}
\revised{\smonth{7} \syear{2010}}

%
\begin{abstract}
We perform an error analysis for numerical approximation methods of
continuous time Markov chain models commonly found in the chemistry
and biochemistry literature. The motivation for the analysis is to
be able to compare the accuracy of different approximation methods
and, specifically, Euler tau-leaping and midpoint tau-leaping. We
perform our analysis under a scaling in which the size of the time
discretization is inversely proportional to some (bounded) power of
the norm of the state of the system. We argue that this is a~more
appropriate scaling than that found in previous error analyses in
which the size of the time discretization goes to zero independent
of the rest of the model. Under the present scaling, we show that
midpoint tau-leaping achieves a higher order of accuracy, in both a~%
weak and a strong sense, than Euler tau-leaping; a result that is in
contrast to previous analyses. We present examples that demonstrate
our findings.
\end{abstract}

%

\end{frontmatter}

%


\section{Introduction}
\label{sec:intro}

This paper provides an error analysis for numerical approximation
methods for continuous time Markov chain models that are becoming
increasingly common in the chemistry and biochemistry literature. Our
goals of the paper are two-fold. First, we want to demonstrate the
importance of considering appropriate scalings in which to carry out
error analyses for the methods of interest. Second, we wish to
provide such an error analysis in order to compare the accuracy of two
different approximation methods. We perform our analysis on the Euler
tau-leaping method first presented in \cite{Gill2001} and a midpoint
tau-leaping method developed below, which is only a slight variant of
one presented in \cite{Gill2001}. The midpoint tau-leaping method
will be demonstrated to be more accurate than Euler tau-leaping in
both a strong and a weak sense, a result that is in contrast to
previous error analyses. We will discuss why previous error analyses
made differing predictions than does ours and argue that the scaling
provided here, or variants thereof, is a more natural and appropriate
choice for error analyses of such methods. We also provide examples
that demonstrate our findings.

\subsection{The basic model}

The motivation for the class of mathematical models under
consideration comes from chemistry and biochemistry, and more
generally from population processes (though we choose the language of
chemistry throughout the paper). We assume the existence of a
chemical reaction system consisting of (i) $d$ chemical species
$\{S_1,S_2,\ldots, S_d\}$ and (ii) a finite set of possible
reactions, which we index by $k$. Each reaction requires some number
of the species as inputs and provides some number of the species as
outputs. For example, the reaction $S_1 \to2S_2$ would require one
molecule of $S_1$ for the input and provide two molecules of $S_2$ for
the output. If reaction $k$ occurs at time $t$, then the state of the
system $X(t)\in\Z^d_{\ge0}$ is updated via addition of the
\textit{reaction vector} $\nu_k \in\Z^d$, which represents the net
change in the abundances of the underlying species:
\[
X(t) = X(t-) + \nu_k.
\]
Returning briefly to the example $S_1 \to2 S_2$, the associated
reaction vector for this reaction would be $[-1,2,0,\ldots,0]^T$.
Finally, we denote by $\nu_k^s$ the vector in $\Z^d_{\ge0}$
representing the source of the $k$th reaction. Returning again to the
example $S_1 \to2S_2$, the source vector for this reaction is
$\nu_k^s = [1,0,\ldots,0]^T$.

We assume that the waiting times for the $k$ reactions are
exponentially distributed with intensity functions $\lambda_k\dvtx
\R^d_{\ge0} \to\R_{\ge0}$. We extend each $\lambda_k$ to all of
$\R^d$ by setting it to zero outside $\R^d_{\ge0}$. This model is a
continuous time Markov chain in $\Z^d_{\ge0}$ with generator
%
%
\begin{equation}\label{eq:generator_exact}
(\EG f)(x) = \sum_{k} \lambda_k(x)\bigl(f(x + \nu_k) - f(x)\bigr),
\end{equation}
where $f \dvtx \Z^d \to\R$ is arbitrary. Kolmogorov's forward equation
for this model, termed the ``chemical master equation'' in the
chemistry and biology literature, is
\[
\frac{d}{dt} P(x, t | \pi) = \sum_k P(x-
\nu_k,t|\pi)\lambda_k(x-\nu_k) - \sum_k P(x,t|\pi)\lambda_k(x),
\]
where\vspace*{1pt} for $x \in\Z^d_{\ge0}$ $P(x,t|\pi)$ represents the probability
that $X(t) = x$, conditioned upon the initial distribution $\pi$. One
representation for path-wise solutions to this model uses a random
time change of Poisson processes,
%
%
\begin{equation}\label{eq:RTC_exact}
X(t) = X(0) + \sum_k Y_k \biggl( \int_0^t \lambda_k(X(s))\,ds
\biggr)\nu_k,
\end{equation}
where the $Y_k$ are independent, unit-rate Poisson processes (see,
e.g., \cite{KurtzPop81}). Note that $\tilde{X}(t) \eqdef X(t) -
\sum_k \int_0^t \lambda_k(X(s)) \,ds\, \nu_k $ is\vspace*{1pt} a martingale with
quadratic covariation matrix $[X]_t = \sum_k Y_k ( \int_0^t
\lambda_k(X(s))\,ds )\nu_k \nu_k^T$.

A common choice of intensity function for chemical reaction systems,
and the one we adopt throughout, is mass action kinetics. Under
mass\vadjust{\goodbreak}
action kinetics, the intensity function for the $k$th reaction is
%
%
\begin{eqnarray}\label{eq:stoch_MA}
\lambda_k(x) &=& \tilde{c}_k \pmatrix{x\cr x-\nu_k^s} =
\frac{\tilde{c}_k}{\prod_{\ell=1}^d \nu_{k \ell}^s!} \prod_{\ell=
1}^d \frac{x_{\ell}!}{(x_{\ell} - \nu^s_{k \ell})!}
1_{\{x_{\ell} \ge0 \}} \nonumber\\[-8pt]\\[-8pt]
&\eqdef& c_k \prod_{\ell= 1}^d
\frac{x_{\ell}!}{(x_{\ell} - \nu^s_{k \ell})!} 1_{\{x_{\ell}
\ge
0 \}},\nonumber
\end{eqnarray}
where $\tilde{c}_k$ is a positive constant and $c_k$ is defined by the
above equation. Mass action kinetics arises by thinking of
$\tilde{c}_k \Delta t$ as the approximate probability that a
\textit{particular} set of the molecules needed in the $k$th reaction
will react over a time-period of size $\Delta t$, and then counting
the number of ways such a reaction could happen. Implicit in the
assumption of mass action kinetics is that the vessel under
consideration is ``well stirred.'' For ease of notation, we will
henceforth drop the indicator functions from our representation of
mass action kinetics. More general rates will be discussed in the
remark at the top of page six.


\subsection{Numerical methods}

There are a number of numerical methods that produce statistically
exact sample paths for the model described above. These include the
stochastic simulation algorithm, better known as Gillespie's algorithm
\cite{Gill76,Gill77}, the first reaction method \cite{Gill76} and
the next reaction method \cite{Anderson2007a,gibson2000}. All such
algorithms perform the same two basic steps multiple times until a
sample path is produced over a desired time interval: first,
conditioned on the current state of the system the amount of time that
passes until the next reaction takes place, $\Delta$, is computed and
second the specific reaction that has taken place is found. If,
however, $\sum_k \lambda_k(X(t)) \gg0$ then $\Delta\approx(\sum_k
\lambda_k(X(t)))^{-1} \ll1$ and the time needed to produce a single
exact sample path over a time interval can be prohibitive.

The approximate algorithm ``tau-leaping'' was developed by Dan
Gillespie in \cite{Gill2001} in an effort to overcome the problem that
$\Delta$ may be prohibitively small. The basic idea of tau-leaping is
to hold the intensity functions fixed over the time interval $[t_n,
t_n+h]$ at the values $\lambda_k(X(t_n))$, where $X(t_n)$ is the
current state of the system, and, under this assumption, compute the
number of times each reaction takes place over this period. As the
waiting times for the reactions are exponentially distributed, this
leads to the following algorithm.
\begin{algorithm}[(Euler tau-leaping)] Set $Z(0) = X(0)$, $t_0 = 0$,
$n=0$ and repeat the following until $t_{n+1} > T$.
\begin{longlist}[(1)]
\item[(1)] Set $Z(t_{n+1}) = Z(t_{n}) + \sum_k
\mathcal{P}_{k,n}(\lambda_k(Z(t_{n})) h) \nu_k$, set $t_{n+1} =
t_{n} + h$ and set $n = n+1$, where $\mathcal{P}_{k,n}(x)$ are
independent Poisson random variables with parameters $x$.
\end{longlist}
\end{algorithm}

Several improvements and modifications have been made to the basic
algorithm described above over the years. However, they are mainly
concerned with how to choose the step-size adaptively
\cite{Cao2006,Gill2003} and/or how to ensure that population values do
not go
negative during the course of a simulation
\cite{Anderson2007b,Cao2005,Chatterjee2005}, and are not explicitly
relevant to the
current discussion.

Similar to (\ref{eq:RTC_exact}), a path-wise representation of Euler
tau-leaping
can be given through a random time change of Poisson processes,
%
%
\begin{equation}\label{eq:RTC_tau}
Z(t) = X(0) + \sum_k Y_k \biggl( \int_0^t \lambda_k\bigl(Z \circ\eta(s)\bigr)
\,ds \biggr)\nu_k,
\end{equation}
where $\eta(s) = t_n$ if $t_n \le s < t_{n+1}$ and the $Y_k$ are as
before. Noting that $\int_0^{t_{n+1}} \lambda_k(Z\circ\eta(s)) \,ds =
\sum_{i=0}^n \lambda_k(Z(t_i))(t_{i+1} - t_i)$ explains our choice to
call this method ``Euler tau-leaping.'' Defining the operator
%
%
\begin{equation}\label{eq:tau_gen}
(\TG_z f)(x) = \sum_k \lambda_k(z)\bigl(f(x + \nu_k) - f(x)\bigr),
\end{equation}
we see that for $t>0$
%
%
\begin{equation}\label{eq:tau_martingale}
\E f(Z(t)) = \E f\bigl(Z\circ\eta(t)\bigr) + \E\int_{\eta(t)}^t \bigl(\TG
_{Z\circ
\eta(t)} f\bigr)(Z(s)) \,ds,
\end{equation}
so long as the expectations exist. Further, we note that
$\tilde{Z}(t) \eqdef Z(t) - \sum_k \int_0^t \lambda_k(Z \circ\eta(s))
\,ds \nu_k$ is a martingale with quadratic covariation matrix
\[
[\tilde{Z}]_t
= \sum_k Y_k\biggl( \int_0^t \lambda_k \bigl(Z\circ\eta(s)\bigr) \,ds
\biggr)\nu_k
\nu_k^T.
\]

It is natural to believe that a midpoint type method would be more
accurate than an Euler type method in many situations. We therefore
define the function
\[
\rho(z) \eqdef z + \frac{1}{2}h\sum_k \lambda_k(z) \nu_k,
\]
which computes an approximate midpoint for the system assuming the
state of the system is $z$ and the time-step is $h$.
\begin{algorithm}[(Midpoint tau-leaping)] Set $\MT(0) = X(0)$, $t_0 =
0$, $n=0$ and repeat the following until $t_{n+1} > T$.
\begin{longlist}[(1)]
\item[(1)] Set $\MT(t_{n+1}) = \MT(t_{n}) + \sum_k
\mathcal{P}_{k,n}(\lambda_k \circ\rho\circ\MT(t_{n}) h) \nu_k$,
set $t_{n+1} = t_{n} + h$ and set $n = n+1$, where
$\mathcal{P}_{k,n}(x)$ are independent Poisson random variables
with parameters $x$.
\end{longlist}
\end{algorithm}

Similar to (\ref{eq:RTC_exact}) and (\ref{eq:RTC_tau}), $\MT(t)$ can
be represented via a random time change of Poisson processes:
\[
\MT(t) = X(0) + \sum_k Y_k \biggl( \int_0^t \lambda_k \circ\rho
\circ\MT\circ\eta(s) \,ds \biggr)\nu_k,
\]
where $\eta(\cdot)$ is as before. For $\TG_z$ defined via
(\ref{eq:tau_gen}) and any $0<t$ and any function $f$
%
%
\begin{equation}\label{eq:mdpt_martingale}
\E f(\MT(t)) = \E f\bigl(\MT\circ\eta(t)\bigr) + \E\int_{\eta(t)}^t
\bigl(\TG_{\rho\circ\MT\circ\eta(t)} f\bigr)(\MT(s)) \,ds.
\end{equation}
Finally, $\tilde{\MT}(t) \eqdef\MT(t) - \sum_k \int_0^t \lambda_k
\circ\rho\circ\MT\circ\eta(s) \,ds\, \nu_k $ is a martingale with
quadratic covariation matrix $[\tilde{\MT}]_t = \sum_k Y_k(
\int_0^t
\lambda_k \circ\rho\circ\MT\circ\eta(s) \,ds )\nu_k
\nu_k^T$. The main goal of this paper is to show that the midpoint
tau-leaping algorithm is indeed more accurate than the Euler
tau-leaping method under an appropriate, and natural, scaling
described in Section \ref{sec:assump_scales}.
%

\begin{remark*}
Historically, the time discretization parameter for tau-leaping has
been $\tau$, thus giving the method its name. We choose to break
from this tradition and denote our time-step by $h$ so as not to
confuse $\tau$ with a~stopping time.
\end{remark*}


\subsection{Previous error analyses}

Under the scaling $h \to0$, Rathinam et al.~\cite{Rathinam2005} performed a consistency
check of Euler tau-leaping and found that the local truncation error
was $\mathcal{O}(h^2)$ for all moments. They also
showed that under this same scaling Euler tau-leaping is first order
accurate in a weak sense in the special case that the intensity
functions $\lambda_k$ are linear \cite{Rathinam2005}. Li extended
these results by showing that as $h \to0$ Euler tau-leaping has a~%
strong error (in the $L^2$ norm) of order $1/2$ and a weak error of
order one~\cite{Li2007}, which agree with classical results pertaining
to numerical analysis of SDEs driven by Brownian motions (see, e.g.,
\cite{KloedenPlaten92}).

Under the scaling $h \to0$, it is readily seen that midpoint
tau-leaping is no more accurate than Euler tau-leaping. This follows
since midpoint tau-leaping consists of making an $\mathcal{O}(h^2)$
correction to the intensity functions used in Euler tau-leaping. As
$h\to0$, this correction becomes negligible as Poisson processes
``ignore'' $\mathcal{O}(h^2)$ corrections, and the accuracy of the two
methods will be the same.

We simply note that while the analyses performed in
\cite{Rathinam2005} and \cite{Li2007} and the argument made in the
previous paragraph are technically correct, performing an analysis as
$h \to0$, independent of the rest of the model, is at odds with the
useful regime of tau-leaping. That is, tau-leaping would only be used
in a regime where $h \gg\Delta$, where $\Delta$ is the expected
amount of time between reactions, for otherwise an exact method would
be performed. Therefore, we should require that
%
%
\begin{equation}\label{eq:h_requirement}
h\gg\biggl(\sum_k \lambda_k(Z(t))\biggr)^{-1} \quad\mbox{or}\quad h \sum_k
\lambda_k(Z(t)) \gg1,
\end{equation}
where $Z(t)$ is the state of the system. In Section
\ref{sec:assump_scales}, we will present a natural scaling for the
models under consideration that does satisfy (\ref{eq:h_requirement})
and under which we will perform our analysis.

\subsection{Paper outline}

The remainder of the paper is organized as follows. In Section
\ref{sec:assump_scales}, we give some natural assumptions on the models
considered in this paper and introduce the scaling under which we
perform our analysis. In Section~\ref{sec:strong}, we perform a strong
error analysis for both the Euler and midpoint tau-leaping methods and
show that midpoint tau-leaping is the more accurate of the two under
our scaling. In Section \ref{sec:weak}, we perform a weak error
analysis of the different methods and again conclude that the midpoint
method is more accurate. In Section \ref{sec:examples}, we present
numerical examples demonstrating our results.


\section{Assumptions on the model}
\label{sec:assump_scales}


\subsection{Scalings of the model and the algorithms}
\label{sec:scaling}

As discussed in the \hyperref[sec:intro]{Intro-} \hyperref[sec:intro]{duction}, tau-leaping methods will only be of
use if the time-discretization parameter $h$ satisfies $h \sum_k
\lambda_k(Z(t)) \gg1$ while $(\sum_k \lambda_k(Z(t)))^{-1} \ll1$,
where~$Z(t)$ is the state of the system at time $t$. There are a
number of ways for the second condition to hold and a modeling choice
must be made.
We make the following natural assumptions:
\begin{longlist}
\item The initial abundance of each species scales with $V$ for some
$V \gg1$.
%
\item Each rate constant satisfies $c_k^V = \mathcal{O}(V^{1 -
\nu_k^s\cdot\vec1})$, where $\vec1 = [1,1,\ldots,1]^T$. In
particular, $c_k^V = d_k/V^{1 - \nu_k^s\cdot\vec1}$ for some
$d_k>0$.
\end{longlist}
We will denote by $X^V$ the normalized process defined as the vector
of abundances divided by $V$, and will denote by $\lambda_k^V$ the
intensity function defined to be mass action kinetics with rate
constants $c_k^V$. This scaling is the so called ``classical
scaling'' and arises naturally by thinking of $V$ as the volume of the
vessel in which the reactions are taking place multiplied by
Avogadro's number \cite{Kurtz72}. In this case, $X^V$~gives the
concentration of each species in moles per unit volume. To understand
the scaling for the rate constants, consider the case of a reaction
requiring as input two constituent molecules: $S_1$ and $S_2$.
Perhaps $S_1 + S_2 \to S_3$. It is reasonable to assume that the
probability that a \textit{particular} pair of $S_1$ and $S_2$
molecules meet, and react, in a~small time interval is inversely
proportional to the volume of the vessel. This same type of logic
holds for the cases in which more than two molecules are needed for a
reaction to take place (i.e., the probability that three particular
molecules meet and react is inversely proportional to the volume
squared). For the case that only one molecule is needed for a
reaction to take place, it is reasonable to assume that the
probability of such a reaction taking place in the next small interval
of time for a \textit{particular} molecule should not scale with the
volume. See also~\cite{Wilkinson2006}, Chapter 6.

Models that satisfy assumptions (i) and (ii) above have an
important property that we will detail here and make use of later.
Let $x(t)$ denote the solution to the deterministic initial value
problem
%
%
\begin{equation}\label{eq:ode_approx}
\dot x(t) = F(x(t)) \eqdef\sum_k d_k x(t)^{\nu_k^s} \nu_k,\qquad x(0)
= x_0 \in\R^d_{\ge0},
\end{equation}
where $d_k$ is defined in assumption (ii) above, and where for any
two vectors $u^v \eqdef u_1^{v_1}\cdots u_d^{v_d}$ and we adopt the
convention that $0^0 = 1$. That is, $x(t)$ is the solution to the
corresponding deterministically modeled chemical reaction system with
mass action kinetics. It was shown in \cite{Kurtz72,Kurtz78} that
for any $\varepsilon> 0$ and any $T>0$, if $X^V(0) = x(0) = x_0$, then
%
%
\begin{equation}\label{eq:Kurtz_limit}
\lim_{V \to\infty} P\Bigl\{\sup_{t \in[0,T]}|X^V(t) - x(t)| \ge
\varepsilon\Bigr\} \to0.
\end{equation}

Denoting $\lambda_k$ as
\textit{deterministic} mass action kinetics with rate constant $d_k$,
it is an exercise to check that for any reaction, that is, zeroth order,
first order, second order, etc., and any $x \in\R^d_{\ge0}$
\[
\lambda_k^{V}(Vx) =
V \lambda_k(x) +
\zeta_k^V(x),
\]
where $\zeta_k^V$ is uniformly bounded in $V$ and is nonzero only if
the reaction requires more than one molecule of a particular species
as an input. For example, for the second order reaction $S_1 + S_2
\to S_3$ we have
\[
\lambda_k^{V}(Vx) = \frac{d_k}{V}
(Vx_1 ) (V x_2) = V d_k
x_1 x_2 = V
\lambda_k(x),
\]
whereas for the second order reaction $2S_1 \to S_3$ we have
\[
\lambda_k^V(Vx) = \frac{d_k}{V}V
x_1( V x_1 - 1) = V d_kx_1^2 -
d_k x_1 = V \lambda_k(x) + \zeta_k^V(x)
\]
with $\zeta_k^V(x) = - d_kx_1$. The term $\zeta_k^V$ will have a true
$V$ dependence if three or more molecules of a particular species are
required as input. We now state the definition $A_k^V(x) \eqdef
\frac{1}{V} \lambda_k^V(Vx)$, and note that for all $x \in\R^d_{\ge
0}$
%
%
\begin{equation}\label{eq:AkDef}
A_k^V(x) = \lambda_k(x) + \frac{1}{V}\zeta_k^V(x)
\end{equation}
and $A_k^V(x) \equiv0$ if $x \notin\R^d_{\ge0}$. Manipulating the
definition of $A_k^V$ shows that for all $x \in\R^d$
%
%
\begin{equation}\label{eq:nec_scaling}
\lambda_k^V(Vx) = V A^V_k(x).
\end{equation}
%

\begin{remark*}
The assumption of mass action kinetics is not critical to the
analysis carried out in this paper. Instead, what is critical to
this particular analysis is that our kinetics satisfies the scaling
(\ref{eq:nec_scaling}) for $A_k^V$ satisfying (\ref{eq:AkDef}) with
$\lambda_k$ sufficiently smooth.
\end{remark*}

We now choose a discretization parameter for the approximate methods
that is dependent upon the assumptions of the model set out above. We
let
%
%
\begin{equation}\label{eq:tau_scale}
h^V \eqdef1/V^{\beta},
\end{equation}
where $0 < \beta< 1$. We note that this scaling satisfies the
necessary requirements detailed above as
\begin{eqnarray*}
\biggl(\sum_k \lambda_k^V(Vx) \biggr)^{-1} &=& \mathcal{O}(V^{-1})
\ll1,\\
V^{-\beta} \biggl( \sum_k \lambda_k^V(Vx)\biggr) &=&
\mathcal{O}(V^{1 - \beta}) \gg1.
\end{eqnarray*}

With this choice of time-step, we let $Z^V$ and $\MT^V$ denote the
processes generated by Euler and midpoint tau-leaping, respectively,
normalized by $V$. We can now state more clearly what the analysis of
this paper will entail. We will consider the case of $V \gg0$ by
letting $V \to\infty$ and consider the relationship of the normalized
approximate processes $Z^V$ and $\MT^V$ to the original process $X^V$,
normalized similarly. Note that all three processes converge to the
solution of (\ref{eq:ode_approx}). We will perform both weak and
strong error analyses. In the strong error analysis, we will consider
$L^1$ convergence as opposed to the more standard (at least for
systems driven by Brownian motions) $L^2$ convergence. The reason for
this is simple: the It\^o isometry makes working with the $L^2$-norm
easier in the Brownian motion case, whereas Poisson processes lend
themselves naturally to analysis in the $L^1$-norm.

We remark that it is clear that the choice of scaling laid out in this
section and assumed throughout the paper will not explicitly cover all
cases of interest. For example, one may choose to use approximation
methods when (i) the abundances of only a strict subset of the
constituent species are in an~$\mathcal{O}(V)$ scaling regime, or (ii)
it is the rate constants themselves that are $\mathcal{O}(V)$ while the
abundances are $\mathcal{O}(1)$, or (iii) there is a mixture of the
previous two cases with potentially more than two natural scales in
the system. Our analysis will not be directly applicable to such
cases. However, the purpose of this analysis is not to handle every
conceivable case. Instead, our purpose is to try and give a more
accurate picture of how different tau-leaping methods approximate the
exact solution, both strongly and weakly, in at least one plausible
setting and we believe that the analysis detailed in this paper
achieves this aim. Further, we believe that error analyses conducted
under different modeling assumptions can be carried out in similar
fashion.


\subsection{Redefining the kinetics}
\label{sec:red_kin}

Before proceeding to the analysis, we allow ourselves one change to
the model detailed in the previous section. As we will be considering
approximation methods in which changes to the state of the system are
determined by Poisson random variables (which can produce arbitrarily
large values), there will always be a positive probability that
solutions will leave a region in which the scaling detailed above is
valid. Multiple options are available to handle such a situation.
One option would be to define a stopping time for when the process
leaves a predetermined region in which the scaling regime is valid and
then only perform the analysis up to that stopping time. Another
option, and the one we choose, is to simply modify the kinetics by
multiplying by a cutoff function that makes the intensity functions
zero outside such a region. This has the added benefit of
guaranteeing the existence of all moments of the processes involved.
Note that without this truncation or some other additional assumption
guaranteeing the existence of the necessary moments, some of the
moment estimates that follow may fail; however, the convergence in
probability and convergence in distribution results in Theorems
\ref{thm:exact_asympt_euler} and \ref{thm:mdpt_exact} would still be
valid.

Let $\gamma\ge0$ be $C^{\infty}$ with compact support
$\Omega_{\gamma} \subset\R^d_{>0}$, with $\gamma(x) = 1$ for all $x
\in B_r(x(t))$ for some $r>0$, where $x(t)$ satisfies
(\ref{eq:ode_approx}). Now, we redefine our intensity functions by
setting
%
%
\begin{equation} \label{eq:red_kin}
\lambda_k^V(x) = \gamma(x/V)c_k^V \prod_{\ell= 1}^d
\frac{x_{\ell}!}{(x_{\ell} - \nu_{k \ell}^s)!}\qquad \mbox{for }
x \in
\R^d,
\end{equation}
where $c_k^V$ still satisfies the scaling detailed in the previous
section. It is easy to check that the redefined kinetics still
satisfies $\lambda_k^V(Vx) = VA_k^V(x)$, where now $A_k^V(x)$ has also
been redefined by multiplication by $\gamma(x)$. Further, the
redefined $\lambda_k^V$ is identical to the previous function on the
domain of interest to us. That is, they only differ if the process
leaves the scaling regime of interest. For the remainder of the paper,
we assume our intensity functions are given by (\ref{eq:red_kin}).
Finally, we note that for each $k$ we have the existence of an $L_k>0$
such that
%
%
\begin{equation}\label{eq:Lk_bound}
\sup_{x \in\R^d, |\alpha| < \infty} |D^{\alpha} A_k^V(x)| \le L_k.
\end{equation}


\section{Strong error analysis for Euler and midpoint tau-leaping}
\label{sec:strong}

Throughout this section, we assume a time discretization
$0=t_0 < t_1 <\cdots< t_N = T$ with $t_n - t_{n-1} = h^V = V^{-\beta
}$ for
some $0 < \beta< 1$. In Section \ref{sec:preliminaries} we give some
necessary technical results. In Section \ref{sec:bounds} we give
bounds for $\sup_{t \le T} \E|X^V(t) - Z^V(t)|$ and $\sup_{t \le T}
\E|X^V(t) - \MT^V(t)|$ in terms of $V$, where $X^V(t), Z^V(t)$ and~%
$\MT^V(t)$ are the normalized processes and satisfy the representations
%
%
\begin{eqnarray}
\label{eq:exact_scaled}
X^V(t) &=& X^V(0) + \frac{1}{V}\sum_k Y_k\biggl(V \int_0^t
A_k^V(X^V(s)) \,ds \biggr)\nu_k,\\
\label{eq:tau_scaled}
Z^V(t) &=& X^V(0) + \frac{1}{V}\sum_k Y_k\biggl(V \int_0^t A_k^V\bigl(Z^V
\circ\eta(s)\bigr) \,ds \biggr)\nu_k, \\
\label{eq:mdpt_scaled}
\MT^V(t) &=& X^V(0) + \frac{1}{V}\sum_k Y_k\biggl(V \int_0^t A_k^V
\circ\rho^V \circ\MT^{V} \circ\eta(s) \,ds
\biggr)\nu_k,
\end{eqnarray}
where
\[
\rho^V(z) \eqdef z + \frac{1}{2}V^{-\beta}\sum_k A_k^V(z) \nu_k
\]
and $\eta(s) = t_n$ for $s \in[t_n,t_{n+1})$.
In Sections \ref{sec:exact_euler} and \ref{sec:exact_mdpt}, we use
different couplings of the processes than those above to provide the
exact asymptotics of the error processes $X^V - Z^V$ and $X^V - \MT^V$.


\subsection{Preliminaries}
\label{sec:preliminaries}

We present some technical, preliminary concepts that will be used
ubiquitously throughout the section. For a more thorough reference of
the material presented here, see \cite{Kurtz86}, Chapter 6. We begin
by defining the following filtrations that are generated by the
Poisson processes $Y_k$:
\begin{eqnarray*}
\mathcal{F}_{\tilde u} &\eqdef& \sigma\{Y_k(s_k) \dvtx s_k \le u_k\},\\
\mathcal{F}_u^i &\eqdef& \sigma\{Y_{k}(r), Y_{i}(s) \dvtx k \ne i,
s \le u, r \in[0,\infty) \},
\end{eqnarray*}
where $\tilde{u}$ is a multi-index and $u$ is a scalar.

\begin{lemma}\label{lem:stopping_time}
Suppose that $X(t)$ satisfies (\ref{eq:RTC_exact}) with nonnegative
intensity functions $\lambda_k$. For $t \ge0$ and a choice of $k$,
%
%
\begin{equation}\label{eq:tau_k}
\tau_k(t) = \int_0^t \lambda_k(X(s)) \,ds
\end{equation}
is an $\{\mathcal{F}_u^k\}$-stopping time.
\end{lemma}
\begin{pf}
For $u \ge0$, let $\alpha(u)$ satisfy
\[
\int_0^{\alpha(u)} \lambda_{k}(X(s))\,ds = u,
\]
where we take\vspace*{2pt} $\alpha(u) = \infty$ if $\int_0^{\infty}
\lambda_k(X(s))\,ds < u$. Then $\alpha(u)$ is adapted to~%
$\mathcal{F}_u^k$ and $\{\tau_k(t) \le u\} = \{t \le\alpha(u)\} \in
\mathcal{F}_u^k$.
\end{pf}

Therefore, if the processes $X(t)$ and $Z(t)$ satisfy
(\ref{eq:RTC_exact}) with nonnegative intensity functions
$\lambda_{k,1}$ and $\lambda_{k,2}$, respectively, then for $t,s \ge
0$ and a choice of $k$,
%
%
\begin{eqnarray} \label{eq:stopping}
&&\E\biggl| Y_k\biggl( \int_0^t \lambda_{k,1}(X(r))\,dr \biggr) -
Y_k\biggl( \int_0^s \lambda_{k,2}(Z(r))\,dr \biggr) \biggr| \nonumber\\[-8pt]\\[-8pt]
&&\qquad= \E
\biggl| \int_0^t \lambda_{k,1}(X(r))\,dr - \int_0^s
\lambda_{k,2}(Z(r)) \,dr \biggr|,\nonumber
\end{eqnarray}
because (i) both the maximum and minimum of two stopping times are
stopping times, and (ii) $Y_k$ is monotone.

Similarly to above, one can show that $\tau(t) \eqdef(\tau_1(t),
\tau_2(t),\ldots)$, where $\tau_k(t)$ is as in (\ref{eq:tau_k}), is a
multi-parameter $\{\mathcal{F}_{\tilde u}\}$-stopping time. We now
define the filtration
\[
\mathcal{G}_t \eqdef\mathcal{F}_{\tau(t)}
\]
and note that by the conditions of Section \ref{sec:red_kin} the
centered process
%
%
\begin{equation}\label{eq:center}
\qquad\tilde{Y}_k\biggl(\int_0^t \lambda_k(X(s))\,ds \biggr) \eqdef
Y_k\biggl(\int_0^t \lambda_k(X(s))\,ds \biggr) - \int_0^t \lambda_k(X(s))\,ds
\end{equation}
is a square integrable martingale, with respect to $\mathcal{G}_t$,
with quadratic variation $ Y_k(\int_0^t
\lambda_k(X(s))\,ds ) $. This fact will be used repeatedly
throughout the paper.


\subsection{Bounds on the strong error}
\label{sec:bounds}
The following theorems give bounds on the errors $\sup_{t \le T}
\E|X^V(t) - Z^V(t)|$ and $\sup_{t \le T} \E|X^V(t) - \MT^V(t)|$.
\begin{theorem}\label{thm:strong_tau_bound}
Let $X^V(t)$ and $Z^V(t)$ satisfy (\ref{eq:exact_scaled}) and
(\ref{eq:tau_scaled}), respectively, for $t \le T$. Then there
exists a constant $C = C(T) > 0$ such that
\[
\sup_{t \le T} \E|X^V(t) - Z^V(t)| \le CV^{-\beta} = Ch^V.
\]
\end{theorem}

\begin{pf} For $t \in[0,T]$, define $E(t)\eqdef\E|X^V(t) -
Z^V(t)|$. Using (\ref{eq:stopping})\break and~(\ref{eq:Lk_bound}),
\begin{eqnarray*}
E(t) &\le& \biggl(\sum_k |\nu_k| L_k\biggr) \E\int_0^t | X^V(s) -
Z^V \circ\eta(s) | \,ds \\
&\le& \biggl( \sum_k|\nu_k| L_k\biggr) \int_0^t E(s) \,ds + \biggl(\sum_k |\nu_k| L_k\biggr)
\E\int_0^t | Z^V(s) - Z^V \circ\eta(s) | \,ds.
\end{eqnarray*}
The second term on the right above can be bounded similarly,
\[
\E\int_0^t | Z^V(s) - Z^V \circ\eta(s) | \,ds \le
\biggl( \sum_k |\nu_k| L_k\biggr) t V^{-\beta},
\]
and the result holds via Gronwall's inequality.\vspace*{-2pt}
\end{pf}
%
%
\begin{theorem}\label{thm:strong_mdpt_bound}
Let $X^V(t)$ and $\MT^V(t)$ satisfy (\ref{eq:exact_scaled}) and
(\ref{eq:mdpt_scaled}), respectively, for $t \le T$. Then there
exists a constant $C = C(T) > 0$ such that
\[
\sup_{t \le T} \E|X^V(t) - \MT^V(t)| \le
CV^{-\kappa(\beta)}\qquad \mbox{where }
\kappa(\beta) = \min\biggl\{\frac{1 + \beta}{2}, 2\beta\biggr\}.\vspace*{-2pt}
\]
\end{theorem}

Before proving Theorem \ref{thm:strong_mdpt_bound}, we present some
preliminary material. Let $F^V(x) \eqdef\sum_k A_k^V(x) \nu_k$ and define
\[
U^{V,1}(s)\eqdef\MT^V(s) -\rho^V\circ\MT^V\circ\eta(s)=\MT^
V(s)-\MT^V\circ\eta(s)-\tfrac{1}{2} V^{-\beta} F^V\bigl(\MT^V \circ
\eta(s)\bigr)
\]
and
\[
\tilde{U}^{V,1}(s) \eqdef\bigl(s-\eta(s)- \tfrac{1}{2} V^{-\beta
}\bigr)F^V\bigl(\MT
^V\circ
\eta(s)\bigr).
\]
Then
%
%
\begin{eqnarray}\label{eq:udiff}
&&
U^{V,1}(s)-\tilde{U}^{V,1}(s)\nonumber\\
&&\qquad = \tilde{\MT}^V(s)-\tilde{\MT
}^V\circ
\eta
(s)\\
&&\qquad\quad{} + \bigl(s-\eta(s)\bigr)\bigl(F^V\bigl(\rho^V\circ\MT^V\circ
\eta(s)\bigr)-F^V\bigl(\MT^V\circ\eta(s)\bigr)\bigr),\nonumber
\end{eqnarray}
where $\tilde{\MT}^V(t) \eqdef\MT^V(t)-\int_0^tF^V(\rho^V\circ
\MT^V\circ\eta(s))\,ds$ is a martingale.\vspace*{-2pt}
\begin{lemma}\label{lem:udif_one}
For all $0<\beta<1$, there exists a $C>0$ such that
\[
\sup_{s \le\infty} \E|U^{V,1}(s) - \tilde{U}^{V,1}(s)| \le
CV^{-\kappa
(\beta)}.\vspace*{-2pt}
\]
\end{lemma}
\begin{pf}
Clearly, the third term on the right-hand side of (\ref{eq:udiff}) is $\mathcal{
O}(V^{-2\beta})$ uniformly in $s$. Thus,
\begin{eqnarray*}
\E|U^{V,1}(s) - \tilde{U}^{V,1}(s)| &\le& \E|\tilde\MT^V(s)-\tilde
{\MT}^V\circ\eta(s)| +c_1V^{-2\beta}\\
&\le& \biggl( \frac{1}{V} \sum_k|\nu_k|^2\E\int_{\eta
(s)}^sA_k^V\bigl(\rho^V\circ
\MT^V\circ\eta(r)\bigr) \,dr \biggr)^{1/2}\\
&&{} + c_1V^{-2\beta}\\
&\le& c_2 V^{-(1+\beta)/2} + c_1 V^{-2\beta}
\end{eqnarray*}
for constants $c_1$ and $c_2$ which do not depend upon $s$.\vspace*{-2pt}
\end{pf}
\begin{lemma}\label{lem:udifest}
For all $0<\beta<1$ and $0 < t$, and for $\alpha\in\{2,3,4,\ldots\}$
\[
\lim_{V\rightarrow\infty}V^{\alpha\beta}\sup_{s\leq t}\E
[|U^{V,1}(s)-\tilde
{U}^{V,1}(s)|^{\alpha}]=0.\vspace*{-2pt}\vadjust{\goodbreak}
\]
\end{lemma}
\begin{pf}
The third term on the right-hand side of (\ref{eq:udiff}) is
$\mathcal{O}(V^{-2\beta})$, so
\begin{eqnarray*}
\E|U^{V,1}(s)-\tilde{U}^{V,1}(s)|^2&\leq&C\bigl(\E|\tilde
\MT^V(s)-\tilde{\MT}^V\circ\eta(s)|^2+V^{-4\beta}\bigr)\\
&\leq&\frac{C}{V} \sum_k|\nu_k|^2\E\int_{\eta
(s)}^sA_k^V\bigl(\rho^V\circ
\MT^V\circ\eta(r)\bigr) \,dr + CV^{-4\beta}\\
&=&\mathcal{O}\bigl(V^{-((1+\beta) \wedge4\beta)}\bigr)
\end{eqnarray*}
showing the $\alpha= 2$ case.

It is simple to show that $V^{\alpha\beta}\sup_{s \le t}\E[
|U^{V,1}(s)-\tilde{U}^{V,1}(s)|^{\alpha}]$ is uniformly bounded in
$V$ for
any $\alpha\in\Z_{\ge0}$. The $\alpha= 2$ case then gives the
necessary bounds for the arbitrary $\alpha$ case.
\end{pf}

Note that by Lemmas \ref{lem:udif_one} and \ref{lem:udifest}
%
%
\begin{eqnarray} \label{eq:decomp}
&&
A_k^V(\MT^V(s))-A_k^V\bigl(\rho^V\circ\MT^V\circ\eta(s)\bigr) \nonumber\\
&&\qquad= \nabla
A_k^V\bigl(\rho^V\circ\MT^V\circ\eta(s)\bigr) \cdot U^{V,1}(s) + \mathcal{O}
(V^{-2\beta})\\
&&\qquad= \nabla A_k^V\bigl(\rho^V\circ\MT^V\circ\eta(s)\bigr)\cdot\tilde{U}^{V,1}(s)
+ \mathcal{O}\bigl(V^{-\kappa(\beta)}\bigr).\nonumber
\end{eqnarray}
We finally note that for any bounded function $g$ and any $n\ge0$
\[
\int_{t_n}^{t_{n+1}} g(\eta(s))\tilde{U}^{V,1}(s) \,ds = 0
\]
and so for any $t > 0$
%
%
\begin{eqnarray} \label{eq:final_bound_U}
&&\int_0^t g(\eta(s))\tilde{U}^{V,1}(s) \,ds\nonumber\\
&&\qquad= \frac{1}{8}\bigl( \bigl(2t -
2\eta(t) - V^{-\beta}\bigr)^2 - V^{-2\beta}\bigr) g(\eta(t))F^V\bigl(\MT^V
\circ\eta(t)\bigr)\\
&&\qquad= \mathcal{O}(V^{-2\beta}).\nonumber
\end{eqnarray}
\begin{pf*}{Proof of Theorem \ref{thm:strong_mdpt_bound}} For $t \le T$
define $E(t) \eqdef\E|X^V(t) - \MT^V(t)|$. Letting $c_i$ denote
constants
\begin{eqnarray*}
E(t) &\le& \sum_k |\nu_k| \E\biggl| \int_0^t A_k^V(X^V(s)) \,ds -
\int_0^t
A_k^V \circ\rho^V \bigl(\MT^{V} \circ\eta(s)\bigr) \,ds \biggr| \\
&\le& c_1 \int_0^t E(s) \,ds + \sum_k |\nu_k| \E\biggl| \int_0^t
A_k^V(\MT^V(s)) - A_k^V \circ\rho^V \bigl(\MT^{V} \circ\eta(s)\bigr) \,ds
\biggr| \\
&\le& c_1 \int_0^t E(s) \,ds + c_2 V^{-\kappa(\beta)},
\end{eqnarray*}
where the final inequality used both (\ref{eq:decomp}) and
(\ref{eq:final_bound_U}). The result now follows from Gronwall's
inequality.
\end{pf*}


\subsection{Exact asymptotics for Euler tau-leaping}
\label{sec:exact_euler}

Throughout this section and the next, all convergences are understood
to hold on bounded intervals. More explicitly, we write $X^V \to X$
if $\lim_{V \to\infty} P\{\sup_{t \le T}|X^V(t) - X(t)|>\varepsilon\} =
0$ for all $\varepsilon> 0$ and $T>0$. Because of the simplifying
assumptions made on the kinetics in Section \ref{sec:red_kin}, it is
not difficult to show that $X^V \to X$ also implies $\lim_{V\to\infty}
\E\sup_{t \le T} |X^V(t) - X(t)| = 0$. In light of this, when we
write $X^V = Z^V + \mathcal{O}(V^{-p})$ for some $p>0$ in this section
and the next we mean that for any $T>0$ there exists a $C(T)$ such
that
\[
\lim_{V \to\infty} V^{p}\E\sup_{t \le T} |X^V(t) - Z^V(t)|
\le C(T).
\]
Finally, recall that $F^V(x) = \sum_k A_k^V(x) \nu_k$ and note that
the function $F(x)$ and the deterministic process $x(s)$ used in the
characterization of the error processes are defined via
(\ref{eq:ode_approx}).

Theorem \ref{thm:strong_tau_bound} suggests that $X^V - Z^V$ scales
like $V^{-\beta}$. In this section, we make this precise by
characterizing the limiting behavior of $V^{\beta} (X^V - Z^V)$, as \mbox{$V
\to\infty$}. To get the exact asymptotics for the Euler tau-leap
method, we will use the following coupling of the processes involved:
%
%
\begin{eqnarray}
\label{eq:exactcouple_tau}
\hspace*{22pt}X^V(t) &=& X^V(0)\nonumber\\
\hspace*{22pt}&&{}+ \frac{1}{V} \sum_k \biggl[Y_{k,1} \biggl(
V\int_0^tA_k^V(X^V(s)) \wedge A_k^V\bigl(Z^V\circ\eta(s)\bigr)
\,ds\biggr)\nonumber\\[-8pt]\\[-8pt]
\hspace*{22pt}&&\hspace*{72.8pt}\hspace*{-30.3pt}{} +Y_{k,2} \biggl( V\int_0^t
A_k^V(X^V(s))\nonumber\\
\hspace*{22pt}&&\hspace*{176.7pt}\hspace*{-46.5pt}\hspace*{-30.3pt}{}
-A_k^V(X^V(s))\wedge A_k^V\bigl(Z^V \circ\eta(s)\bigr)\,ds \biggr) \biggr]
\nu_k,
\nonumber
\\
\label{eq:tau_couple}
\hspace*{22pt}Z^V(t)&=&X^V(0)\nonumber\\
\hspace*{22pt}&&{}+\frac{1}{V} \sum_k \biggl[ Y_{k,1} \biggl( V\int_0^t
A_k^V(X^V(s))\wedge A_k^V\bigl(Z^V\circ\eta(s)\bigr) \,ds\biggr)
\nonumber\\[-8pt]\\[-8pt]
\hspace*{22pt}&&\hspace*{72.8pt}\hspace*{-30.3pt}{} +Y_{k,3}\biggl(V\int_0^tA_k^V\bigl(Z^V\circ
\eta(s)\bigr)\nonumber\\
\hspace*{22pt}&&\hspace*{176.7pt}\hspace*{-75.3pt}{}-A_k^V(X^V(s))\wedge A_k^V\bigl(Z^V \circ\eta
(s)\bigr)\,ds\biggr)\biggr]\nu_k.
\nonumber
\end{eqnarray}
It is important to note that the distributions of $X^V$ and $Z^V$
defined via~(\ref{eq:exactcouple_tau}) and (\ref{eq:tau_couple}) are
the same as those for the processes defined via
(\ref{eq:exact_scaled}) and (\ref{eq:tau_scaled}).

The following lemma is easy to prove using Doob's inequality.
\begin{lemma}\label{lem:obvious_lemma1}
For $X^V$ and $Z^V$ given by (\ref{eq:exactcouple_tau}) and
(\ref{eq:tau_couple}), $X^V - Z^V \to0$.
\end{lemma}

Combining Lemma \ref{lem:obvious_lemma1} and (\ref{eq:Kurtz_limit})
shows that $Z^V - x \to0$, where $x$ is the solution to the
associated ODE. Similarly, $Z^V \circ\eta- x \to0$. These facts
will be used throughout this section.

Centering the Poisson processes, we have
%
%
\begin{eqnarray}\label{eq:errpr}
X^V(t)-Z^V(t)&=&M^V(t)+\int^t_0 F^V(X^V(s))-F^V\bigl(Z^V\circ\eta(s)\bigr)
\,ds \nonumber\\
&=& M^V(t)+\int^t_0 F^V(X^V(s))-F^V(Z^V(s)) \,ds\\
&&{} +\int^t_0
F^V(Z^V(s))-F^V\bigl(Z^V\circ\eta(s)\bigr)\,ds,\nonumber
\end{eqnarray}
where $M^V$ is a martingale.

To obtain the desired results, we must understand the behavior of the
first and third terms on the right-hand side of (\ref{eq:errpr}). We begin\vspace*{1pt} by
considering the third term. We begin by defining
$U^{V,2}$ and $\tilde{U}^{V,2}$ by
\[
U^{V,2}(s) \eqdef Z^V(s) - Z^V \circ\eta(s),\qquad \tilde{U}^{V,2}(s)
\eqdef\bigl(s - \eta(s)\bigr)F^V\bigl(Z^V \circ\eta(s)\bigr).
\]
Then,
\[
U^{V,2}(s) - \tilde{U}^{V,2}(s) = \tilde{Z}^V(s) - \tilde{Z}^V \circ
\eta(s),
\]
where $\tilde{Z}^V(t) \eqdef Z^V(t) - \int_0^tF^V(Z^V \circ\eta(s))
\,ds$ is a martingale. Thus,
%
%
\begin{eqnarray}\label{eq:tau_F_expand}\quad
&&
F^V(Z^V(s)) - F^V\bigl(Z^V \circ\eta(s)\bigr)\nonumber\\
&&\qquad= DF^V\bigl(Z^V \circ
\eta(s)\bigr)U^{V,2}(s) + \mathcal{O}(V^{-2\beta})
\nonumber\\[-8pt]\\[-8pt]
&&\qquad= DF^V\bigl(Z^V \circ\eta(s)\bigr)\tilde{U}^{V,2}(s) + DF^V\bigl(Z^V \circ
\eta(s)\bigr)\bigl(U^{V,2}(s) - \tilde{U}^{V,2}(s)\bigr)\nonumber\\
&&\qquad\quad{} + \mathcal{
O}(V^{-2\beta}).\nonumber
\end{eqnarray}

\begin{lemma}\label{lem:tau_U}
For all $0<\beta<1$, $0<t$, and $\alpha\in\{2,3,4,\ldots\}$
\[
\lim_{V \to\infty} V^{\alpha\beta} \sup_{s\le t} \E[|U^{V,2}(s) -
\tilde{U}^{V,2}(s)|^{\alpha}] = 0.
\]
\end{lemma}
\begin{pf}
The proof is similar to that of Lemma \ref{lem:udifest}.
\end{pf}

We may now characterize the limiting behavior of the third term of
(\ref{eq:errpr}).

\begin{lemma}\label{lem:tau_U2}
For $0<\beta<1$ and any $t>0$,
\[
V^{\beta} \int_0^t F^V(Z^V(s)) - F^V\bigl(Z^V \circ\eta(s)\bigr) \,ds \to
\frac{1}{2} \int_0^t DF(x(s))F(x(s)) \,ds.
\]
\end{lemma}
\begin{pf}
By (\ref{eq:tau_F_expand}) and Lemma \ref{lem:tau_U}
\begin{eqnarray*}
&&V^{\beta} \int_0^t F^V(Z^V(s)) - F^V\bigl(Z^V \circ\eta(s)\bigr) \,ds\\[-2pt]
&&\qquad =
V^{\beta}\int_0^t DF^V\bigl(Z^V \circ\eta(s)\bigr)F^V\bigl(Z^V \circ\eta(s)\bigr)\bigl(s -
\eta(s)\bigr)\,ds + \varepsilon_1^V(t),
\end{eqnarray*}
where $\varepsilon_1^V \to0$ as $V\to\infty$. By Lemma
\ref{lem:obvious_lemma1} convergence results similar to
(\ref{eq:Kurtz_limit}) hold for the process $Z^V \circ\eta$, and
because $\int_{\eta(s)}^{\eta(s) + V^{-\beta}} (r - \eta(s))\,dr =
\frac{1}{2}V^{-2\beta}$, the lemma holds as stated.\vspace*{-2pt}
\end{pf}

Turning now to $M^V$, we observe that the quadratic covariation is
\[
[M^V]_t = \frac{1}{V^2} \sum_k \bigl(N^V_{k,2}(t) + N_{k,3}^V(t)\bigr) \nu_k
\nu_k^T,
\]
where
\begin{eqnarray*}
N^V_{k,2}(t) &\eqdef& Y_k\biggl(V \int_0^t A_k^V(X^V(s)) -
A_k^V(X^V(s))\wedge A_k^V\bigl(Z^V \circ\eta(s)\bigr)\biggr),\\[-2pt]
N^V_{k,3}(t) &\eqdef& Y_k\biggl(V \int_0^t A_k^V\bigl(Z^V \circ\eta(s)\bigr) -
A_k^V(X^V(s))\wedge A_k^V\bigl(Z^V \circ\eta(s)\bigr)\biggr),
\end{eqnarray*}
which as $V \to\infty$ is asymptotic to
%
%
\begin{equation}\label{eq:asympt}
\frac{1}{V}\sum_k \int_0^t \bigl|A_k^V(X^V(s)) - A_k^V\bigl(Z^V \circ\eta(s)\bigr)
\bigr| \,ds\, \nu_k \nu_k^T.
\end{equation}
We have the following lemma.\vspace*{-2pt}
\begin{lemma}\label{lem:euler_mart}
For $0<\beta<1$, $V^{\beta}M^V \to0$, as $V \to\infty$.\vspace*{-2pt}
\end{lemma}
\begin{pf}
Multiplying (\ref{eq:errpr}) by $V^{\alpha}$, we see that
$V^{\alpha} (X^V- Z^V) \rightarrow0$ provided $\alpha<\beta$ (so
that the third term on the right goes to zero) and provided
$V^{\alpha}M^V\rightarrow0$. By the martingale central limit
theorem, the latter convergence holds provided
$V^{2\alpha}[M^V]\rightarrow0$ (see Lemma \ref{lem:MCLT} in
the \hyperref[sec:app]{Appendix}). Let $\alpha_0=\sup\{\alpha\dvtx \alpha\le
\beta, V^{ 2\alpha}[M^V]\rightarrow0\}$. Since $\alpha_0<1$, we
have that $2\alpha_0-1<\alpha_0 \le\beta$, which implies by the
definition of $\alpha_0$ that $V^{2\alpha_0 - 1}(X^V - Z^V) \to0$.
Therefore,
\begin{eqnarray*}
&&
V^{2\alpha_0}[M^V]_t \\[-2pt]
&&\qquad\approx \sum_k \int_0^t V^{2\alpha_0 - 1}
\bigl|A_k^V(X^V(s)) - A_k^V\bigl(Z^V \circ\eta(s)\bigr) \bigr| \,ds\, \nu_k \nu_k^T\\[-2pt]
&&\qquad\approx \sum_k \int_0^t V^{2\alpha_0 - 1} \bigl|\nabla
A_k^V\bigl(Z^V\circ\eta(s)\bigr) \cdot\bigl(Z^V(s) - Z^V \circ\eta(s)\bigr)\bigr|
\,ds\,
\nu_k \nu_k^T\\[-2pt]
&&\qquad\approx \sum_k \int_0^t V^{2\alpha_0 - 1} \bigl|\nabla A_k^V\bigl(Z^V\circ
\eta(s)\bigr) \cdot F^V\bigl(Z^V\circ\eta(s)\bigr)\bigr| \bigl(s - \eta(s)\bigr) \,ds\, \nu_k
\nu_k^T,
\end{eqnarray*}
where in the second approximation we used that $V^{2\alpha_0 -
1}(X^V - Z^V) \to0$, in the third approximation we substituted
$\tilde{U}^{V,2}(s)$ for $U^{V,2}(s)$, and by $f \approx g$ we mean $f
- g
\to0$ as $V \to\infty$. The last expression goes to zero whenever
$2\alpha_0-1 < \beta$, hence the convergence holds.\vspace*{-2pt}
\end{pf}

We now have the following theorem characterizing the behavior of\break
\mbox{$V^{\beta}(X^V\,{-}\,Z^V)$}.\vspace*{-2pt}
\begin{theorem}\label{thm:exact_asympt_euler}
For $X^V$ and $Z^V$ given by (\ref{eq:exactcouple_tau}) and
(\ref{eq:tau_couple}) and for $0<\beta<1$, $V^{\beta}(X^V - Z^V)
\to\Err$, where $\Err$ is the solution to
%
%
\begin{eqnarray}\label{eq:det}
\Err(t) &=& \int_0^t DF(x(s))\Err(s)\,ds\nonumber\\[-10pt]\\[-10pt]
&&{} + \frac{1}{2}\int_0^t
DF(x(s))F(x(s))\,ds,\qquad \Err(0) = 0.\nonumber\vspace*{-2pt}
\end{eqnarray}
\end{theorem}
\begin{pf}
Multiply (\ref{eq:errpr}) by $V^{\beta}$ and observe that
\[
V^{\beta} \int_0^t F^V(X^V(s)) - F^V(Z^V(s)) \,ds \approx\int_0^t
DF^V(Z^V(s))V^{\beta}\bigl(X^V(s) - Z^V(s)\bigr) \,ds.
\]
The theorem now follows directly from Lemmas \ref{lem:tau_U2} and
\ref{lem:euler_mart}.\vspace*{-2pt}
\end{pf}

%




\subsection{Exact asymptotics for midpoint tau-leaping}
\label{sec:exact_mdpt}

Throughout this section, the Hessian matrix associated with a real
valued function $g$ will be denoted by $Hg$. Also, for any vector
$U$, we will denote by $U^THF^V(x)U$ the vector whose $i$th component
is $U^THF_i^VU$, and similarly for $F$.

The goal of this section is to characterize the limiting behavior
of
\[
V^{\kappa(\beta)}\bigl(X^V(t)-{\MT}^V(t)\bigr),
\]
where
\[
\kappa(\beta)= \min\{2\beta,(1+\beta)/2\}=\cases{
2\beta, &\quad $\beta<1/3$,\cr
(1+\beta)/2, &\quad $\beta\geq1/3$.}
\]
To get the exact asymptotics for the midpoint method, we will use the
following representation of the processes involved:
%
%
\begin{eqnarray}\label{exactcouple}
X^V(t) &=& X^V(0)\nonumber\\[-2pt]
&&{}+ \frac{1}{V} \sum_k \biggl[Y_{k,1}
\biggl(V\int_0^tA_k^V(X^ V(s)) \wedge A_k^V\bigl(\rho^V\circ
\MT^V\circ\eta(s)\bigr)\,ds\biggr)\nonumber\\[-10pt]\\[-10pt]
&&\hphantom{{}+ \frac{1}{V} \sum_k \biggl[}{}
 +Y_{k,2}\biggl(V\int_0^tA_k^V(X^V(s))\nonumber\\[-1pt]
&&\hspace*{104.2pt}{}-A_k^V(X^V(s))\wedge
A_k^V\bigl(\rho^V\circ\MT^V\circ\eta(s)\bigr)\,ds\biggr)\biggr]\nu
_k,\nonumber\\[-1pt]
\label{mdcouple}
{\MT}^V(t)&=&X^V(0)\nonumber\\
&&{}+\frac{1}{V} \sum_k \biggl[ Y_{k,1}
\biggl(V\int_0^t A_k^V(X^V(s))\wedge A_k^V\bigl(\rho^V\circ\MT^V\circ\eta
(s)\bigr)\,ds\biggr)\nonumber\\[-8pt]\\[-8pt]
&&\hphantom{{}+ \frac{1}{V} \sum_k \biggl[}{}
+Y_{k,3}\biggl(V\int_0^tA_k^V\bigl(\rho^V\circ\MT^V\circ\eta
(s)\bigr)\nonumber\\
&&\hspace*{104.2pt}{}-A_k^V(X^V(s))\wedge A_k^V\bigl(\rho^V\circ\MT^V\circ\eta
(s)\bigr)\,ds\biggr)\biggr]\nu_k.\nonumber
\end{eqnarray}

The following is similar to Lemma \ref{lem:obvious_lemma1}.
\begin{lemma}\label{lem:obvious_lemma2}
$\!$For $X^V$ and $\MT^V$ given by (\ref{exactcouple}) and
(\ref{mdcouple}), \mbox{$X^V - \MT^V \to0$}.
\end{lemma}

Combining Lemma \ref{lem:obvious_lemma2} and (\ref{eq:Kurtz_limit})
shows that $\MT^V - x \to0$, where $x$ is the solution to the
associated ODE. Similarly $\MT^V \circ\eta- x \to0$. These facts
will be used throughout this section.

Centering the Poisson processes, we have
%
%
\begin{eqnarray}\label{errpr}\quad
X^V(t)-\MT^V(t)&=& M^V(t)+\int^t_0 F^V(X^V(s))-F^V\bigl(\rho^V\circ
\MT^V\circ\eta(s)\bigr) \,ds \nonumber\\
&=& M^V(t)+\int^t_0 F^V(X^V(s))-F^V(\MT^V(s)) \,ds\\
&&{} +\int^t_0
F^V(\MT^V(s))-F^V\bigl(\rho^V\circ\MT^V\circ\eta(s)\bigr)\,ds,\nonumber
\end{eqnarray}
where $M^V$ is a martingale.

As before, we must understand the behavior of the first and third
terms on the right-hand side of (\ref{errpr}). We begin by considering the
third term. Proceeding as in the previous sections, we define $U^{V,3}$
and $\tilde U^{V,3}$ as
\[
U^{V,3}(s)\eqdef\MT^V(s) -\rho^V\circ\MT^V\circ\eta(s)=\MT^
V(s)-\MT^V\circ\eta(s)-\tfrac{1}{2} V^{-\beta} F^V\bigl(\MT^V \circ
\eta(s)\bigr)
\]
and
\[
\tilde{U}^{V,3}(s) \eqdef\bigl(s-\eta(s)- \tfrac{1}{2} V^{-\beta
}\bigr)F^V\bigl(\MT
^V\circ
\eta(s)\bigr).
\]
Then
%
%
\begin{eqnarray}\label{udiff}
&&
U^{V,3}(s)-\tilde{U}^{V,3}(s)\nonumber\\[2pt]
&&\qquad = \tilde{\MT}^V(s)-\tilde{\MT
}^V\circ
\eta
(s)\\[2pt]
&&\qquad\quad{} + \bigl(s-\eta(s)\bigr)\bigl(F^V\bigl(\rho^V\circ\MT^V\circ
\eta(s)\bigr)-F^V\bigl(\MT^V\circ\eta(s)\bigr)\bigr),\nonumber
\end{eqnarray}
where $\tilde{\MT}^V(t) \eqdef\MT^V(t)-\int_0^tF^V(\rho^V\circ
\MT^V\circ\eta(s))\,ds$ is a martingale. Then
%
%
\begin{eqnarray} \label{decomp}
&&
F^V(\MT^V(s))-F^V\bigl(\rho^V\circ\MT^V\circ\eta(s)\bigr) \nonumber\\[2pt]
&&\qquad= DF^
V\bigl(\rho^V\circ\MT^V\circ\eta(s)\bigr)U^{V,3}(s) \nonumber\\[2pt]
&&\qquad\quad{} +\tfrac{1}{2} U^{V,3}(s)^THF^V\bigl(\rho^V\circ
\MT^V\circ\eta(s)\bigr) U^{V,3}(s) +
\mathcal{O}(V^{-3\beta})\nonumber\\[-7pt]\\[-7pt]
&&\qquad= DF^V\bigl(\rho^V\circ\MT^V\circ\eta(s)\bigr)\tilde{U}^{V,3}(s) \nonumber\\[2pt]
&&\qquad\quad{} + \tfrac{1}{2}
U^{V,3}(s)^THF^V\bigl(\rho^V\circ\mathcal{Z}^V\circ\eta
(s)\bigr)U^{V,3}(s) \nonumber\\[2pt]
&&\qquad\quad{} +DF^V\bigl(\rho^V\circ\MT^V\circ\eta
(s)\bigr)\bigl(U^{V,3}(s)-\tilde
{ U}^{V,3}(s)\bigr)+\mathcal{O}(V^{-3\beta}).\nonumber
\end{eqnarray}

\begin{lemma}\label{udifest}
For all $0<\beta<1$, $0 < t$, and $\alpha\in\{2,3,4,\ldots\}$
\[
\lim_{V\rightarrow\infty}V^{\alpha\beta}\sup_{s\leq t}\E
[|U^{V,3}(s)-\tilde
{U}^{V,3}(s)|^{\alpha}]=0.
\]
\end{lemma}
\begin{pf} The proof is similar to Lemma \ref{lem:udifest}.
\end{pf}

Let
\[
\kappa_1(\beta)= \min\{2\beta,\beta+1/2\}=\cases{
2\beta, &\quad $\beta<1/2$,\cr
\beta+1/2, &\quad $\beta\geq1/2$.}
\]
Note that $\kappa_1(\beta) \ge\kappa(\beta)$ for all $\beta\ge0$.
\begin{lemma}
For $0<\beta<\frac{1}{2}$ and each $t>0$,
%
%
\begin{eqnarray}\label{asm1}
&&V^{2\beta}\int_0^tDF^V\bigl(\rho^V\circ\MT^V\circ
\eta(s)\bigr)\bigl(U^{V,3}(s)-\tilde{U}^{V,3}(s)\bigr) \,ds\nonumber\\[-8pt]\\[-8pt]
&&\qquad\rightarrow\frac{1}{4}
\int_0^t DF(x(s))^2F(x(s))\,ds,\nonumber
\end{eqnarray}
for $\beta=\frac{1}{2}$
%
%
\begin{eqnarray}\label{asm2}
&&V\int_0^tDF^V\bigl(\rho^V\circ\MT^V\circ\eta(s)
\bigr)\bigl(U^{V,3}(s)-\tilde{U}^{V,3}(s)\bigr)\,ds\nonumber\\[-8pt]\\[-8pt]
&&\qquad\Rightarrow M_1(t) + \frac{1}{4}
\int_0^tDF(x(s))^2F(x(s))\,ds,\nonumber
\end{eqnarray}
and for $\frac{1}{2}<\beta<1$,
%
%
\begin{equation}\label{asm3}
\quad\qquad
V^{\beta+ 1/2}\int_0^tDF^V\bigl(\rho^V\circ\MT^
V\circ\eta(s)\bigr)\bigl(U^{V,3}(s)-\tilde{U}^{V,3}(s)\bigr)\,ds\Rightarrow
M_1(t),
\end{equation}
where $M_1$ is a mean zero Gaussian process with independent
increments and quadratic covariation
%
%
\begin{equation}\label{eq:quad_var}
[M_1]_t=\frac{1}{3} \int_0^t \sum_kA_k(x(s))DF(x(s)) \nu_k \nu_k^T
DF(x(s))^Tds.
\end{equation}
\end{lemma}
\begin{pf}
By Lemma \ref{mgincrem} in the \hyperref[sec:app]{Appendix},
\begin{eqnarray*}
M_1^V(t)&\eqdef&\int_0^tDF^V\bigl(\rho^V\circ\MT^V\circ\eta(s)\bigr)\bigl(\tilde
{\MT}^V(s)-\tilde{\MT}^V\circ\eta(s)\bigr)\,ds\\
&&{} +DF^V\bigl(\rho^V\circ\MT^V\circ\eta(t)\bigr)\bigl(\tilde{\MT}^
V(t)-\tilde{\MT}^V\circ\eta(t)\bigr)\bigl(\eta(t)+V^{-\beta}-t\bigr)
\end{eqnarray*}
is a martingale and its quadratic covariation matrix is
\[
\int_0^t\bigl(\eta(s)+V^{-\beta}-s\bigr)^2DF^V\bigl(\rho^V\circ\MT^V\circ\eta
(s)\bigr)\,d[\tilde{\MT}^V]_sDF^V\bigl(\rho^V\circ\MT^V\circ\eta(s)\bigr)^T.
\]
Noting that $\int_{\eta(s)}^{\eta(s)+V^{-\beta}}(\eta(r)+V^{-\beta}
-r)^2\,dr=\frac{1}{3} V^{-3\beta}$, it follows that
\[
V^{2\beta+1}[M_1^V]_t\rightarrow\frac{1}{3}\int_0^t\sum_kA_k(x(s)
)DF(x(s))\nu_k\nu_k^TDF(x(s))^Tds,
\]
so by the martingale central limit theorem $V^{\beta+1/2}M_1^V$
converges in distribution to a mean zero Gaussian process with
independent increments and quadratic variation (\ref{eq:quad_var}).

Since $V^{1/2}(\tilde{\MT}^V-\tilde{\MT}^V\circ\eta) \to0$, the
integral on the left-hand side of (\ref{asm1}), (\ref{asm2}) and (\ref{asm3})
can be replaced by
%
%
\begin{eqnarray}\label{change}
&&M_1^V(t)+ \int_0^tDF^V\bigl(\rho^V\circ\MT^V\circ\eta(s)\bigr)\bigl(s-\eta
(s)\bigr)\nonumber\\[-8pt]\\[-8pt]
&&\qquad\hspace*{33pt}{}\times\bigl(F^V\bigl(\rho^V\circ\MT^V\circ\eta(s)\bigr)-F^V\bigl(\MT^V\circ
\eta(s)\bigr)\bigr) \,ds\nonumber
\end{eqnarray}
without\vspace*{1pt} changing the limits. The second term in (\ref{change})
multiplied by $V^{2\beta}$ converges to $\frac{1}{4}\int_0^t
DF(x(s))^2F(x(s))\,d s$ on bounded\vspace*{1pt} time intervals and the three limits
follow.
\end{pf}
\begin{lemma}
For $0<\beta<1$,
\begin{eqnarray*}
&&V^{2\beta} \frac{1}{2} \int_0^t U^{V,3}(s)^T HF^V\bigl(\rho^V\circ\MT^V
\circ\eta(s)\bigr)U^{V,3}(s)\,ds \\
&&\qquad\to\frac{1}{24} \int_0^tF(x(s))^T HF(x(s))
F(x(s))\,ds.
\end{eqnarray*}
\end{lemma}
\begin{pf}
By Lemma \ref{udifest}, we can replace $U^{V,3}$by $\tilde{U}^{V,3}$.
Observing that $\int_{\eta(s)}^{\eta(s)+V^{-\beta}}(s-\eta
(s)-\frac{1}{2} V^{-\beta})^2 \,ds = \frac{1}{12}V^{-3\beta}$,
\begin{eqnarray*}
&&V^{2\beta}\frac{1}{2} \int_0^t
\biggl(s-\eta(s)-\frac{1}{2}V^{-\beta}\biggr)^2F^V\bigl(\MT^V\circ\eta(s)\bigr)^T
HF^V\bigl(\rho^V \circ\MT^V\circ\eta(s)\bigr)\\[-1pt]
&&\qquad\hspace*{18pt}{}\times F^V\bigl(\MT^V\circ\eta(s)\bigr)\,ds
\end{eqnarray*}
converges as claimed.
\end{pf}

We may now characterize the behavior of the third term of (\ref{errpr}).
\begin{lemma}\label{zmcorr}
Let
\[
R^V(t)=\int^t_0\biggl(s-\eta(s)-\frac{1}{2} V^{-\beta}\biggr)DF^V\bigl(\rho^V\circ
\MT^V\circ\eta(s)\bigr)F^V\bigl(\MT\circ\eta(s)\bigr)\,ds.
\]
Then for $0<\beta<\frac{1}{2}$,
\begin{eqnarray*}
&& V^{2\beta}\biggl(\int^t_0 \bigl(F^V(\MT^V(s))
-F^V\bigl(\rho^V \circ\MT^V \circ\eta(s)\bigr) \bigr)\,ds-R^V(t) \biggr)\\[-2pt]
&&\qquad \to\frac{1}{4} \int_0^t DF(x(s))^2 F(x(s))\,ds + \frac{1}{24}
\int_0^t F(x(s))^THF(x(s))F(x(s))\,ds,
\end{eqnarray*}
for $\beta=\frac{1}{2}$,
\begin{eqnarray*}
&&V\biggl(\int^t_0\bigl(F^V(\MT^V(s))-F^V\bigl(\rho^V\circ\MT^V\circ
\eta
(s)\bigr)\bigr) \,ds-R^V(t) \biggr) \\[-2pt]
&&\qquad \Rightarrow M_1(t) + \frac{1}{4}\int_0^t DF(x(s))^2F(x(s))\,ds
\\[-2pt]
&&\qquad\quad{}+ \frac{1}{24} \int_0^tF(x(s))^THF(x(s))F(x(s))\,ds
\end{eqnarray*}
and for $\frac{1}{2} < \beta< 1$,
\[
V^{\beta+ 1/2}\int^t_0\bigl(F^V(\MT^V(s))-F^V\bigl(\rho^V\circ
\MT^V\circ\eta(s)\bigr)\bigr)\,ds \Rightarrow M_1(t).
\]
\end{lemma}
\begin{remark*}Note that $V^{2\beta}R^V$is uniformly bounded,
$R^V\circ\eta\equiv0$, and
\[
R^V(t) = \tfrac{1}{2}\bigl[\bigl(t - \eta(t)\bigr)^2 - \bigl(t - \eta(t)\bigr)V^{-\beta
}\bigr]
DF^V\bigl(\rho^V\circ\MT^V\circ\eta(t)\bigr)F^V\bigl(\MT\circ\eta(t)\bigr).
\]
\end{remark*}
\begin{pf*}{Proof of Lemma \ref{zmcorr}}
The lemma follows from (\ref{decomp}), the previous lemmas, and by
noting that $\int_0^t DF^V(\rho^V \circ\MT^V \circ\eta(s))
\tilde{U}^{V,3} (s) \,ds = R^V(t)$.
\end{pf*}

We now turn to $M^V$ and observe that
\[
[M^V]_t=\frac{1}{V^2}\sum_k\bigl(N_{k,2}^V(t)+N_{k,3}^V(t)\bigr)\nu_k\nu_k^
T,\vadjust{\goodbreak}
\]
where
\begin{eqnarray*}
N_{k,2}^V(t) &\eqdef& Y_k\biggl( V \int_0^t A_k^V(X^V(s)) -
A_k^V(X^V(s))\wedge A_k^V\bigl(\rho^V \circ\MT^V \circ\eta(s)\bigr) \,ds
\biggr),\\
N_{k,3}^V(t) &\eqdef& Y_k\biggl(V \int_0^t A_k^V\bigl(\rho^V \circ\MT^V
\circ
\eta(s)\bigr) \\
&&\hphantom{Y_k\biggl(}\hspace*{22.6pt}{}
- A_k^V(X^V(s))\wedge A_k^V\bigl(\rho^V \circ\MT^V \circ
\eta(s)\bigr) \,ds\biggr),
\end{eqnarray*}
which as $V \to\infty$ is asymptotic to
\[
\frac{1}{V} \sum_k\int_0^t\bigl|A_k^V(X^V(s))-A_k^V\bigl(\rho^V\circ
\MT^V\circ\eta(s)\bigr)\bigr|\,ds\, \nu_k\nu_k^T.
\]
Consequently, we have the following.
\begin{lemma}\label{poisclt}
For $0 < \beta<1$, $V^{(1+\beta)/2} M^V \Rightarrow M$ where
$M$ is a mean-zero Gaussian process with independent increments and
quadratic covariation
\[
[M]_t = \sum_k \frac{1}{4} \int_0^t|\nabla A_k(x(s))\cdot
F(x(s))|\,ds \,\nu_k\nu_
k^T.
\]
\end{lemma}
\begin{pf}
Multiplying (\ref{errpr}) by $V^{\alpha}$, we see that $V^{\alpha}
(X^V-\MT^V)\rightarrow0$ provided $\alpha<\kappa_1(\beta)$ (so
that the third term on the right goes to zero) and provided
$V^{\alpha}M^V\rightarrow0$. By the martingale central limit
theorem, the latter convergence holds provided
$V^{2\alpha}[M^V]\rightarrow0$. Let $\alpha_0=\sup\{\alpha\dvtx \alpha
\le(\beta+1)/2, V^{ 2\alpha}[M^V]\rightarrow0\}$. We make two
observations. First, because $\alpha_0<1$, we have that
$2\alpha_0-1<\alpha_0$. Second, because $\alpha_0 \le(\beta+
1)/2$, we have that $2\alpha_0 - 1 \le\beta$, and, in particular,
$2\alpha_0 - 1 < \kappa_1(\beta)$ for all $\beta\in(0,1)$.
Combining these observations with the definition of $\alpha_0$ shows
that $V^{2(2\alpha_0-1)}[M^V]_t \to0$ and hence
$V^{2\alpha_0-1}(X^V - \mathcal{Z}^V) \to0$. We now have
\begin{eqnarray*}
V^{2\alpha_0}[M^V]_t&\approx&\sum_k\int_0^tV^{2\alpha_0-1}\bigl|A_k(X^
V(s))-A_k\bigl(\rho^V\circ\MT^V\circ\eta(s)\bigr)\bigr|\,ds \,\nu_k\nu_k^T\\
&\approx&\sum_k\int_0^tV^{2\alpha_0-1}\biggl|s-\eta(s)-\frac
12V^{-\beta} \biggr|\\
&&\hphantom{\sum_k\int_0^t}{}
\times\bigl|\nabla A_k\bigl(\rho^V\circ\MT^V\circ\eta(s)\bigr)\cdot
F^V\bigl(\MT^ V\circ\eta(s)\bigr)\bigr|\,ds \,\nu_k\nu_k^T,
\end{eqnarray*}
where in the second line we used that $V^{2\alpha_0 - 1}(X^V - Z^V)
\to0$, and then
substituted $\tilde{U}^{V,3}(s)$ for $U^{V,3}(s)$.
Since the last expression would go to zero if $2\alpha_0-1$ were
less than $\beta$, we see that $2\alpha_0-1=\beta$, that is,
$\alpha_0=(\beta+1)/2$. Furthermore, observing that $\int_{\eta
(s)}^{\eta(s)+V^{-\beta}}|s-\eta(s)-\frac12V^{-\beta} |\,ds=\frac
14V^{-2\beta}$,\vadjust{\goodbreak} we see that
\[
V^{\beta+1}[M^V]_t=V^{2\alpha_0}[M^V]_t\rightarrow\sum_k
\frac{1}{4} \int_0^t|\nabla A_k(x(s)) \cdot
F(x(s))|\,ds \,\nu_k\nu_k^T,
\]
and the lemma follows by the martingale central
limit theorem.
\end{pf}

Collecting the results, we have the following theorem.
\begin{theorem}\label{thm:mdpt_exact}
Let
\[
\mathcal{H}(t) = \frac{1}{6} \int_0^t DF(x(s))^2
F(x(s))\,ds + \frac{1}{24} \int_0^t F(x(s))^THF(x(s)) F(x(s))\,ds.
\]
For $0<\beta<\frac{1}{3}$, $V^{2\beta}(X^V-\MT^V-R^V)\rightarrow
\mathcal{E}_1$, where $\mathcal{E}_1$ is the solution of
%
%
\begin{equation}\label{err1}
\mathcal{E}_1(t)=\int_0^tDF(x(s))\mathcal{E}_1(s)\,ds+\mathcal{H}
(t), \qquad\Err_1(0) = 0.
\end{equation}
For $\beta=\frac{1}{3}$, $V^{2\beta}(X^V-\MT^V-R^V) \Rightarrow
\mathcal{ E}_2$, where $\mathcal{E}_2$ is the solution of
%
%
\begin{equation}\label{err2}
\mathcal{E}_2(t)=M(t)+\int_0^tDF(x(s))\mathcal{E}_2(s)\,d
s+\mathcal{H}(t), \qquad\Err_2(0) = 0.
\end{equation}
For $\frac{1}{3} <\beta<1$, $V^{(1+\beta)/2}(X^V-\MT^V) \Rightarrow
\mathcal{E}_3$, where $\mathcal{E}_3$ is the solution of
%
%
\begin{equation}\label{err3}
\mathcal{E}_3(t)=M(t)+\int_0^tDF(x(s))\mathcal{E}_3(s)\,d s,\qquad
\Err_3(0) = 0.
\end{equation}
\end{theorem}
\begin{pf}
For $\beta\leq\frac{1}{3}$, $R^V$ is $\mathcal{O}(V^{-2\beta})$.
Subtract $R^V$ from both sides of~(\ref{errpr}) and observe that
\begin{eqnarray*}
&&
\int^t_0 F^V(X^V(s))-F^V(^{}\MT^V(s)) \,ds \\
&&\qquad\approx\int_0^tDF^V(
\MT^V(s))\bigl(X^V(s)-\MT^V(s)-R^V(s)\bigr)\,ds\\
&&\qquad\quad{} +\int_0^tDF^V(\MT^V(s))R^V(s)\,ds.
\end{eqnarray*}
Since
\[
V^{2\beta}\int_0^tDF^V(\MT^V(s))R^V(s)\,ds\rightarrow-\frac{1}{12}
\int_
0^t DF(x(s))^2F(x(s))\,ds,
\]
the first\vspace*{1pt} two parts follow from Lemmas \ref{zmcorr} and
\ref{poisclt}.

For $\beta>\frac{1}{3}$, $(1+\beta)/2<2\beta\wedge\kappa_1(\beta)$,
so $ V^{(1+\beta)/2}R^V\rightarrow0$ and
\[
V^{(1+\beta)/2}\int^t_0 F^V(\MT^V(s))-F^V\bigl(\rho^V\circ\MT^
V\circ\eta(s)\bigr) \,ds \rightarrow0,
\]
and the third part follows by Lemma \ref{poisclt}.
\end{pf}


\section{Weak error analysis}
\label{sec:weak}

As in previous sections, we assume the existence of a time
discretization $0=t_0 < t_1 <\cdots< t_N = T$ with $t_n - t_{n-1} =
V^{-\beta}$ for some $0<\beta<1$. We also recall that $\eta(s) = t_n$
for $t_n \le s < t_{n+1}$ for each $n\le N-1$.

Let $X^V$ be a Markov process with generator
%
%
\begin{equation}\label{eq:gen_exact_scaled}
(\EGV f)(x) = \sum_k V A_k^V(x)\bigl(f(x + \nu_k/V) - f(x)\bigr).
\end{equation}
Defining the operator
%
%
\begin{equation}\label{eq:tau_oper_scaled}
(\TGV_z f)(x) = \sum_k VA_k^V(z)\bigl(f(x + \nu_k/V) - f(x)\bigr),
\end{equation}
we suppose that $Z^V$ and $\MT^V$ are processes that satisfy
%
%
\begin{equation}\label{eq:tau_mart_scaled}
\E f(Z(t)) = \E f\bigl(Z\circ\eta(t)\bigr) + \E\int_{\eta(t)}^t \bigl(\TG
_{Z\circ
\eta(t)} f\bigr)(Z(s)) \,ds
\end{equation}
and
%
%
\begin{equation}\label{eq:mdpt_mart_scaled}
\E f(\MT(t)) = \E f\bigl(\MT\circ\eta(t)\bigr) + \E\int_{\eta(t)}^t
\bigl(\TG_{\rho^V \circ\MT\circ\eta(t)} f\bigr)(\MT(s)) \,ds
\end{equation}
for all $t>0$, respectively.

We begin with the weak error analysis of Euler tau-leaping, which
is immediate in light of Theorem \ref{thm:exact_asympt_euler}.
\begin{theorem}\label{thm:weak_tau_exact}
Let $X^V(t)$ be a Markov process with generator
(\ref{eq:gen_exact_scaled}) and let $Z^V(t)$ be a process that
satisfies (\ref{eq:tau_mart_scaled}) for the operator
(\ref{eq:tau_oper_scaled}). Then, for any continuously differentiable
function $f$ and any $t \le T$,
\[
\lim_{V \to\infty} V^{\beta} \bigl( \E f(X^V(t)) - \E f(Z^V(t))
\bigr) = \Err(t) \cdot\nabla f(x(t)),
\]
where $\Err(t)$ satisfies (\ref{eq:det}).
\end{theorem}
\begin{pf}
Without loss of generality, we may assume that $X^V(t)$ and~$Z^V(t)$
satisfy (\ref{eq:exactcouple_tau}) and (\ref{eq:tau_couple}),
respectively. The proof now follows immediately from a combination
of Taylor's theorem and Theorem \ref{thm:exact_asympt_euler}.
\end{pf}
\begin{remark*}
Because the convergence in Theorem \ref{thm:weak_tau_exact} is to a
constant independent of the step-size of the method, we see that
Richardson extrapolation techniques can be carried out. However, we
have not given bounds on the next order correction, and so cannot
say how much more accurate such techniques would be.
\end{remark*}

We now consider the weak error analysis of the midpoint method.
\begin{theorem}\label{thm:weak_mdpt_bound}
Let $X^V(t)$ be a Markov process with generator
(\ref{eq:gen_exact_scaled}) and let $\MT^V(t)$ be a process that
satisfies (\ref{eq:mdpt_mart_scaled}) for the operator
(\ref{eq:tau_oper_scaled}). Then, for any two times continuously\vadjust{\goodbreak}
differentiable function $f$ with compact support, there exists a
constant $C = C(f,T) > 0$ such that
\[
V^{2\beta}|\E f(X^V(T)) - \E f(\MT^V(T))| \le C.
\]
\end{theorem}

Before proving Theorem \ref{thm:weak_mdpt_bound}, some preliminary
material is needed. Let $\LL\eqdef\{y \dvtx y = x/V, x
\in\Z^d\}$, and for $x \in\LL$ and a given function $f$, let
%
%
\begin{equation}\label{eq:v}
v(t,x) = \E_x f(X^V(t)),
\end{equation}
where $\E_x$ represents the expectation conditioned upon $X^V(0) = x$.
Standard results give that $v(t,x)$ satisfies the following initial
value problem (see, e.g.,~\cite{EK2005}, Proposition 1.5)
%
%
\begin{eqnarray}\label{eq:partial_v}
\partial_t v(t,x) &=& \EGV v(t,x)\nonumber\\
&=& \sum_k VA_k^V(x)\bigl(v(t,x + \nu_k/V) -
v(t,x)\bigr),\\
&&\eqntext{v(0,x) = f(x),\qquad x \in\LL.}
\end{eqnarray}
The above equation can be viewed as a linear system by letting $x$
enumerate over $\LL$ and treating $v(t,x) = v_x(t)$ as functions in
time only. It can even be viewed as finite dimensional because of the
conditions on the intensity functions $A_k^V$. That is, recall that
$A_k^V(x) = 0$ for all $x$ outside the bounded set $\Omega_{\gamma}$
(see Section \ref{sec:red_kin}); thus, for any such $x \notin
\Omega_{\gamma}$, $v(t,x) = v_x(t) \equiv f(x)$, for all $t > 0$.

For concreteness, we now let $M$ denote the number of reactions for the
system under consideration. For $k,\ell\in[1,\ldots,M]$ and $x \in
\LL$, let
%
%
\begin{eqnarray}
\label{eq:D_k}
D_k(t,x) &=& V\bigl(v(t,x+\nu_k/V) - v(t,x)\bigr),\\
\label{eq:D_kl}
D_{k\ell}(t,x) &=& V\bigl(D_k(t,x+\nu_{\ell}/V) - D_k(t,x)\bigr)
\end{eqnarray}
represent approximations to the first and second spatial derivatives
of $v(t,x)$, respectively. For notational ease, we have chosen not to
explicitly note the~$V$ dependence of the functions $v(t,x)$,
$D_k(t,x)$ or $D_{k\ell}(t,x)$.

The following lemma, which should be viewed as giving regularity
conditions for $v(t,x)$ in the $x$ variable, is instrumental in the
proof of Theorem~\ref{thm:weak_mdpt_bound}. The proof is delayed
until the end of the section.
\begin{lemma}\label{lemma:reg}
Let $v(t,x)$, $D_k(t,x)$ and $D_{k\ell}(t,x)$ be given by
(\ref{eq:v}), (\ref{eq:D_k}) and (\ref{eq:D_kl}), respectively, and
let $T > 0$. There exists $K_1 > 0$ and $K_2 > 0$ that do not depend upon
$V$ such that
%
%
\begin{eqnarray}
\label{eq:reg1}
\sup_{t \le T} \sup_{k \le M} \sup_{x \in\LL}
|D_k(t,x)| &\le& K_1,\\
\label{eq:reg2}
\sup_{t \le T} \sup_{k,\ell\le M} \sup_{x \in\LL}
|D_{k\ell}(t,x) | &\le& K_2.
\end{eqnarray}
\end{lemma}


We will also need the following lemma, which gives regularity
conditions for $D_k(t,x)$ in the $t$ variable, and whose proof is also
delayed.\vadjust{\goodbreak}

\begin{lemma}\label{lemma:lemma_t}
Let $D_k(t,x)$ be given by (\ref{eq:D_k}). There exists a $K> 0$
that does not depend upon $V$ such that
\[
\sup_{t\le T} \sup_{k \le M} \sup_{x \in\LL}|D_k(t + h,x) - D_k(t,x)|
\le K h
\]
for all $h > 0$.
\end{lemma}
%
%
\begin{pf*}{Proof of Theorem \ref{thm:weak_mdpt_bound}}
%
Define the function $u(t,x)\dvtx [0,T] \times\LL\to\R$ by
%
%
\begin{equation}\label{eq:u}
u(t,x) \eqdef\E_x f\bigl(X^V(T - t)\bigr),
\end{equation}
and for any $w(t,x)\dvtx\R\times\LL\to\R$ we define the operator
$\mathcal{L}$ by
\begin{eqnarray*}
\mathcal{L} w(t,x) &\eqdef&\partial_t w(t,x) + \EGV w(t,x)\\
&=& \partial_t w(t,x) + \sum_k VA_k^V(x)\bigl(w(t,x + \nu_k/V) -
w(t,x)\bigr).
\end{eqnarray*}
Note that $u(t,x) = v(T-t,x)$, where $v(t,x)$ is given by
(\ref{eq:v}), and so by (\ref{eq:partial_v}) $\mathcal{L}u(t,x) = 0$
for $t \in[0,T]$ and $x \in\LL$. We also define the operator
\begin{eqnarray*}
\mathcal{L}_z w(t,x) &\eqdef&\partial_t w(t,x) + \TGV_z w(t,x)\\
&=& \partial_t w(t,x) + \sum_k VA_k^V(z)\bigl(w(t,x + \nu_k/V) - w(t,x)\bigr),
\end{eqnarray*}
so that by virtue of equation (\ref{eq:mdpt_mart_scaled}), for $t
\le T$ and any differentiable (in $t$) function $w(t,x)$
%
%
\begin{eqnarray}\label{eq:MT_mart}
\E w(t, \MT^V(t)) &=& \E w\bigl(\eta(t),\MT^V\circ\eta(t)\bigr) \nonumber\\[-8pt]\\[-8pt]
&&{}
+\int_{\eta(t)}^{t} \E\mathcal{L}_{\rho^V \circ\MT^V \circ
\eta(t)} w(s,\MT^V (s))\,ds.\nonumber
\end{eqnarray}
Recalling (\ref{eq:u}), we see that
\begin{eqnarray*}
\E u(T,\MT^V(T)) &=& \E f(\MT^V(T)),\\
\E u(T,X^V(T)) &=& \E u(0,X^V(0)) = \E f(X^V(T)).
\end{eqnarray*}
Therefore by (\ref{eq:MT_mart}), and using that $X^V(0) = \MT^V(0)$,
\begin{eqnarray*}
&&
\E f(\MT^V(T)) - \E f(X^V(T)) \\
&&\qquad= \E u(T,\MT^V(T)) - \E
u(0,\MT^V(0)) \\
&&\qquad= \sum_{n = 0}^{N - 1} \E u(t_{n+1},\MT^V(t_{n+1})) - \E u(t_n,
\MT^V(t_n)) \\
&&\qquad= \sum_{n = 0}^{N - 1} \E\int_{t_n}^{t_{n+1}}
\mathcal{L}_{\rho^V \circ\MT^V(t_n)} u(s,
\MT^V(s))\,ds.
\end{eqnarray*}
Because $\mathcal{L}u(t,x) \equiv0$ for $t \le T$ and $x \in\LL$
%
%
\begin{eqnarray} \label{eq:need}
&&\E \int_{t_n}^{t_{n+1}} \mathcal{L}_{\rho^V \circ\MT^V(t_n)}
u(s,\MT^V(s)) \,ds \nonumber\\
&&\qquad= \E\int_{t_n}^{t_{n+1}} \mathcal{L}_{\rho^V
\circ\MT^V(t_n)} u(s,\MT^V(s)) -
\mathcal{L}u(s,\MT^V(s))\,ds \nonumber\\
&&\qquad= \sum_{k} \E\int_{t_n}^{t_{n+1}} V\bigl[A_k^V \bigl(\rho^V \circ
\MT^V(t_n) \bigr) - A_k^V( \MT^V(s) ) \bigr] \\
&&\qquad\quad\hphantom{\sum_{k} \E\int_{t_n}^{t_{n+1}}}{}\times
\bigl( u\bigl(s,\MT^V(s) +
\nu_k/V\bigr) - u(s,\MT^V(s)) \bigr) \,ds \nonumber\\
&&\qquad= \sum_{k} \E\int_{t_n}^{t_{n+1}} \bigl[A_k^V \bigl( \rho^V \circ
\MT^V(t_n)\bigr) - A_k^V( \MT^V(s) ) \bigr] D_k\bigl(T-s,\MT^V(s)\bigr)
\,ds.\hspace*{-22pt}\nonumber
\end{eqnarray}
Thus, it is sufficient to prove that each of the integrals in
(\ref{eq:need}) are $\mathcal{O}(V^{-3\beta})$. By Lemma
\ref{lemma:lemma_t}, each integral term in (\ref{eq:need}) can be
replaced by
%
%
\begin{eqnarray} \label{eq:first_bound}
I_k^V(t_n) &\eqdef&\E\int_{t_n}^{t_{n+1}} \bigl[A_k^V \bigl( \rho^V
\circ\MT^V(t_n)\bigr) - A_k^V( \MT^V(s) ) \bigr]\nonumber\\[-8pt]\\[-8pt]
&&\hphantom{\E\int_{t_n}^{t_{n+1}}}{}
\times D_k\bigl(T-t_n,\MT^V(s)\bigr) \,ds.\nonumber
\end{eqnarray}
The remainder of the proof consists of proving that $I_k^V(t_n) =
\mathcal{O}(V^{-3\beta})$.

Letting $g^V_n(x) \eqdef[A_k^V ( \rho^V \circ\MT^V(t_n))\,{-}\,A_k^V(x)] D_k(T-t_n,x)$
and applying~(\ref{eq:mdpt_mart_scaled}) to the integrand in
(\ref{eq:first_bound}) yields
\begin{eqnarray*}
I_k^V(t_n) &=& \E\int_{t_n}^{t_{n+1}} \bigl[ A_k^V\bigl( \rho^V \circ
\MT^V(t_n)\bigr)-A_k^V(\MT^V(t_n)) \bigr]D_k\bigl(T-t_n,\MT^V(t_n)\bigr) \,ds\\
&&{}
+ \sum_j \E\int_{t_n}^{t_{n+1}} \int_{t_n}^s
VA_j^V\bigl(\rho^V\circ\MT^V(t_n)\bigr) \\
&&\hphantom{{}+ \sum_j \E\int_{t_n}^{t_{n+1}} \int_{t_n}^s}
{}\times\bigl(g^V_n\bigl(\MT^V(r) +
\nu_j/V\bigr) - g^V(\MT^V(r))\bigr) \,dr \,ds.
\end{eqnarray*}
We have
\begin{eqnarray*}
A_k^V\bigl(\rho^V \circ\MT^V(t_n)\bigr) &=& A_k^V(\MT^V(t_n))\\
&&{} + \nabla
A_k^V(\MT^V(t_n)) \cdot\frac{1}{2}V^{-\beta} \sum_j
A_j^V(\MT^V(t_n))\nu_j\\
&&{} + \mathcal{O}(V^{-2\beta}).
\end{eqnarray*}
Thus,
%
%
\begin{eqnarray}\qquad
\label{eq:one}
I_k^V(t_n) &=& \sum_j \frac{1}{2}V^{-\beta} \E\int_{t_n}^{t_{n+1}}
\nabla A_k^V(\MT^V(t_n)) \cdot\nu_j
A_j^V(\MT^V(t_n))\nonumber\\[-8pt]\\[-8pt]
&&\hphantom{\sum_j \frac{1}{2}V^{-\beta} \E\int_{t_n}^{t_{n+1}}}
{}\times
D_k\bigl(T-t_n,\MT^V(t_n)\bigr) \,ds + \mathcal{O}(V^{-3\beta})\nonumber\\
\label{eq:two}
&&{} + \sum_j \E\int_{t_n}^{t_{n+1}} \int_{t_n}^s
VA_j^V\bigl(\rho^V\circ\MT^V(t_n)\bigr) \nonumber\\[-8pt]\\[-8pt]
&&\hphantom{{}+ \sum_j \E\int_{t_n}^{t_{n+1}} \int_{t_n}^s}
{}\times\bigl(g^V_n\bigl(\MT^V(r) + \nu_j/V\bigr) -
g^V(\MT^V(r))\bigr) \,dr \,ds.\nonumber
\end{eqnarray}
After some manipulation, the expected value term of (\ref{eq:two})
becomes
\begin{eqnarray*}
&&\E\int_{t_n}^{t_{n+1}} \int_{t_n}^s VA_j^V\bigl(\rho^V \circ
\MT^V(t_n)\bigr)\bigl[A_k^V(\MT^V(r)) - A_k^V\bigl(\MT^V(r) + \nu_j/V\bigr)
\bigr]\\
&&\hphantom{\E\int_{t_n}^{t_{n+1}} \int_{t_n}^s }
{}\times D_k\bigl(T-t_n,\MT^V(r)+\nu_j/V\bigr) \,dr
\,ds \\
&&\qquad{} + \E\int_{t_n}^{t_{n+1}} \int_{t_n}^s
A_j^V\bigl(\rho^V \circ\MT^V(t_n)\bigr)\bigl[ A_k^V\bigl(\rho^V\circ
\MT^V(t_n)\bigr)- A_k^V(\MT^V(r)) \bigr] \\
&&\qquad\hphantom{{}+\E\int_{t_n}^{t_{n+1}} \int_{t_n}^s}
{}\times D_{kj}\bigl(T-t_n,\MT^V(r)\bigr) \,dr
\,ds.
\end{eqnarray*}
By\vspace*{1pt} Lemma \ref{lemma:reg} the last term above is $\mathcal{
O}(V^{-3\beta
})$. Taylor's theorem and the fact that
$A_j^V(\rho^V \circ\MT^V(t_n)) = A_j^V(\MT^V(t_n) + \mathcal{
O}(V^{-\beta})$ then shows\vadjust{\goodbreak} us that the expected value term of
(\ref{eq:two}) is equal to
%
%
\begin{eqnarray}\label{eq:third}
&&
- \E \int_{t_n}^{t_{n+1}} \int_{t_n}^s A_j^V(\MT^V(t_n))\nabla
A_k^V(\MT^V(r))\nonumber\\
&&\hphantom{- \E \int_{t_n}^{t_{n+1}} \int_{t_n}^s}
{}\times\nu_j D_k\bigl(T-t_n,\MT^V(r)+\nu_j/V\bigr) \,dr \,ds
+\mathcal{O}(V^{-3\beta})\nonumber\\[-8pt]\\[-8pt]
&&\qquad= - \E\int_{t_n}^{t_{n+1}} \int_{t_n}^s
A_j^V(\MT^V(t_n))\nabla A_k^V(\MT^V(r))\cdot\nu_j
D_k\bigl(T-t_n,\MT^V(r)\bigr) \,dr \,ds\nonumber\\
&&\qquad\quad{} + \mathcal{O}(V^{-3\beta}),\nonumber
\end{eqnarray}
where the second equality stems from an application of Lemma
\ref{lemma:reg}.\vspace*{1pt}

By Lemma \ref{lemma:reg}, the function $\phi(x) = A_j^V(\MT^V(t_n))
\nabla A_k^V(x) \cdot\nu_j D_k(T-t_n,x)$ satisfies $\sup_{\ell}
|\phi(x + \nu_{\ell}/V) - \phi(x)| = \mathcal{O}(V^{-1})$. Therefore,
applying (\ref{eq:mdpt_mart_scaled})\break to~(\ref{eq:third}) shows that
(\ref{eq:third}) is equal to
\begin{eqnarray*}
&&- \E\int_{t_n}^{t_{n+1}} \int_{t_n}^s
A_j^V(\MT^V(t_n))\nabla A_k^V(\MT^V(t_n))\cdot\nu_j
D_k\bigl(T-t_n,\MT^V(t_n)\bigr) \,dr \,ds\\
&&\qquad{} + \mathcal{O}(V^{-3\beta}).
\end{eqnarray*}
Noting that the sum over $j$ of the above is the negative of
(\ref{eq:one}) plus an~$\mathcal{O}(V^{-3\beta})$ correction concludes
the proof.
\end{pf*}


Theorem \ref{thm:weak_mdpt_bound} can be strengthened in the case of
$\beta< 1/3$.
\begin{theorem}\label{thm:weak_mdpt_exact}
Let $X^V(t)$ be a process with generator (\ref{eq:gen_exact_scaled})
and let~$\MT^V(t)$ be a process that satisfies
(\ref{eq:mdpt_mart_scaled}) for the operator
(\ref{eq:tau_oper_scaled}). Suppose also that $\beta< 1/3$. Then,
for any continuously differentiable function $f$,
\[
\lim_{V \to\infty} V^{2\beta} \bigl( \E f(X^V(T)) - \E f(\MT^V(T))
\bigr) = \Err_1(T) \cdot\nabla f(x(T)),
\]
where $\Err_1(t)$ satisfies (\ref{err1}).
\end{theorem}
\begin{pf} Noting that $R^V(T) \equiv0$, this is an immediate
consequence of Theorem \ref{thm:mdpt_exact}.
\end{pf}
\begin{remark*}
In Theorem \ref{thm:weak_tau_exact}, we provided an explicit
asymptotic value for the scaled error of Euler tau-leaping in terms
of a solution to a differential equation for all scales,
$0<\beta<1$, of the leap step. However, Theorem
\ref{thm:weak_mdpt_exact} gives a similar result for the midpoint
method only in the case $0<\beta<1/3$. For the case $1/3\le\beta<
1$, Theorem \ref{thm:weak_mdpt_bound} only shows that the error is
asymptotically bounded by a constant. The reason for the
discrepancy in results is because in Section \ref{sec:strong} we
were able to show that the dominant component of the pathwise error
for Euler tau-leaping for all $\beta\in(0,1)$ and for midpoint
tau-leaping for $\beta\in(0,1/3)$ was a term that converged to
a~deterministic process. However, in the case $\beta\ge1/3$ for
midpoint tau-leaping, the dominant term of the error is a nonzero
Gaussian process. We note that this random error process should not
be viewed as ``extra fluctuations,'' as they are present in the
other cases. In these other cases, they are just dominated by the
error that arises from the deterministic ``drift'' or ``bias'' of
the error process. We leave the exact characterization of the weak
error of the midpoint method in the case $\beta\ge1/3$ as an open
problem.
\end{remark*}

We now present the delayed proofs of Lemmas \ref{lemma:reg} and
\ref{lemma:lemma_t}.
\begin{pf*}{Proof of Lemma \ref{lemma:reg}}
Let $C_1 > 0$ be such that
\[
\sup_{x\in\LL} \sup_k |D_k(0,x)| =\sup_{x\in\LL} \sup_k \bigl|V\bigl(f(x +
\nu_k/V) -
f(x)\bigr)\bigr| \le C_1.
\]
Using (\ref{eq:partial_v}), a tedious reordering of terms shows that
$D_k(t,x)$ satisfies
%
%
\begin{eqnarray} \label{eq:firstDifference}
\partial_t D_k(t,x) &=& \sum_j A_j^V(x) V [D_k(t,x + \nu_j/V) -
D_k(t,x)]\nonumber\\[-8pt]\\[-8pt]
&&{} + \sum_j \bigl(A_j^V(x + \nu_k/V) -
A_j^V(x)\bigr)VD_j(t,x+\nu_k/V).
\nonumber
\end{eqnarray}
Similarly to viewing $v(t,x) = v_x(t)$ as a finite-dimensional
linear system, (\ref{eq:firstDifference}) can be viewed as a linear
system for the variables $D_k(t,x) = D_{\{k,x\}}(t)$, for $k \in
[1,\ldots,M]$ and $x \in\LL$. Because $A_j^V(x) \equiv0$ for all $x
\notin\Omega_{\gamma}$, we see that $\partial_t D_k(t,x) \equiv0$
for all $x$ such that $x \notin\Omega_{\gamma}$ and $x+\nu_j/V
\notin\Omega_{\gamma}$ for all $j \in[1,\ldots,M]$. Therefore, the
system (\ref{eq:firstDifference}) can be viewed as finite
dimensional also.

Let $\Gamma_1 = [1,\ldots,M] \times\LL$. We enumerate the system
(\ref{eq:firstDifference}) over $b \in\Gamma_1$. That is, for $b =
\{k,x\} \in\Gamma_1$ we let $D_b(t) = D_k(t,x) = D_{b_1}(t,b_2)$.
After some ordering of the set $\Gamma_1$, we let $\R^{\Gamma_1}$
denote the set of (infinite) vectors, $v$, whose $b$th component
is $v_b \in\R$, and then denote $D(t) \in\R^{\Gamma_1}$ as the
vector whose $b$th component is $D_b(t)$. Next, for each $b =
\{k,x\} \in\Gamma_1$, we let
\[
{S}_b \eqdef\sum_j A_j^V(b_2) = \sum_j A_j^V(x)
\]
and let $r_b, R_b \in\R^{\Gamma_1}$ satisfy
\begin{eqnarray*}
R_b \cdot v &=& \sum_j A_j^V(b_2) v_{\{b_1,b_2 + \nu_j/V\}},\\
r_b \cdot v &=& \sum_j \bigl(A_j^V(b_2 + \nu_{b_1}/V) -
A_j^V(b_2)\bigr)Vv_{\{j,b_2 +\nu_{b_1}/V \}}
\end{eqnarray*}
for all $v \in\R^{\Gamma_1}$. It is readily seen that for any $b$
both $R_b$ and $r_b$ have at most~$M$ nonzero components. Also, by
the regularity conditions on the functions~$A_j^V$'s, the absolute
value of the nonzero terms of $r_b$ are uniformly bounded above by
some $K$, which is independent of $V$. Finally, note that $R_b
\cdot1 = S_b$. Combining the previous few sentences shows that for
any vector $v \in\R^{\Gamma_1}$, we have the two inequalities
%
%
\begin{eqnarray}
\label{eq:normhelp}
|R_b \cdot v| &=& \biggl|\sum_j A^V_{j}(b_2) v_{\{b_1,b_2 +
\nu_j/V\}} \biggr|
\le S_b \|v\|_{\infty},\\
\label{eq:normhelp2}
|r_b \cdot v| &\le& KM \| v\|_{\infty},
\end{eqnarray}
where $\| v\|_{\infty} \eqdef\sup_{b \in\Gamma_1}|v_b|$. We now
write (\ref{eq:firstDifference}) as
\[
D_b'(t) = -V S_b D_b(t) + VR_b \cdot D(t) + r_b \cdot D(t),
\]
and so for each $b \in\Gamma_1$
%
%
\begin{equation} \label{eq:normsquared}\qquad
\frac{d}{dt}D_b(t)^2 = -2 V S_b D_b(t)^2 + 2VD_b(t) R_b \cdot
D(t) + 2D_b(t) r_{b} \cdot D(t).
\end{equation}
Only a finite number of the terms $D_b(t)$ are changing in time and
so there is a $b^{1}$ and a $t_1 \in(0,T]$ for which
$|D_{b^{1}}(t)| = \|D(t)\|_{\infty}$ for $t \in[0,t_1]$. By
(\ref{eq:normhelp}), we have that for this $b^{1}$ and any $t \in
[0,t_1]$
\[
\int_0^t D_{b^{1}}(s) R_{b^{1}} \cdot D(s)\,ds \le\int_0^{t}
S_{b^{1}} |D_{b^{1}} (s) | \|D(s)\|_{\infty} \,ds = \int_0^{t}
S_{b^{1}} D_{b^{1}} (s)^2 \,ds,
\]
which, after integrating (\ref{eq:normsquared}), yields
\begin{eqnarray*}
\|D(t)\|^2_{\infty} &=& D_{b^{1}}(t)^2 \le D_{b^{1}}(0)^2 + 2
\int_0^t D_{b^{1}}(s) r_{b^{1}} \cdot D(s) \,ds \\
&\le&
\|D(0)\|^2_{\infty} + 2KM \int_0^{t} \|D(s)\|_{\infty}^2 \,ds,
\end{eqnarray*}
where the final inequality makes use of (\ref{eq:normhelp2}). An
application of Gronwall's inequality now gives us that for $t \in
[0,t_1]$
\[
\|D(t)\|_{\infty}^2 \le\|D(0)\|^2_{\infty}e^{2KMt}.
\]
To complete the proof, continue this process for $i\ge2$ by
choosing the~$b^{i}$ for which $|D_{b^{i}}(t)|$ is maximal on the
time interval $t_i - t_{i-1}$. We must have $\lim_{i \to\infty}
t_i = T$ because (i) there are a finite number of time varying
$D_b(t)$'s and (ii) each~$D_b(t)$ is differentiable. After taking
square roots, we find $\sup_{t \le T} \|D(t)\|_{\infty} \le
\|D(0)\|_{\infty}e^{KMT}\le C_1e^{KMT}$, which is equivalent to~(\ref{eq:reg1}).

We now turn our attention to showing (\ref{eq:reg2}), which we show
in a similar manner. There is a $C_2 > 0$ such that for all $x \in
{\LL}$ and $k, \ell\in[1,\ldots,M]$,
\begin{eqnarray*}
|D_{k\ell}(0,x)| &=& V^2|f(x + \nu_{\ell}/V + \nu_k/V) - f(x +
\nu_{\ell}/V) - f(x + \nu_k/V) + f(x)| \\
&\le& C_2.
\end{eqnarray*}
Another tedious reordering of terms, which makes use of
(\ref{eq:firstDifference}), shows that~$D_{k\ell}(t,x)$ satisfies
\begin{eqnarray*}
&&
\partial_t D_{k\ell}(t,x) \\
&&\qquad= \sum_j A_j^V(x)V[D_{k\ell}(t,x +
\nu_j/V) - D_{k\ell}(t,x)] \\
&&\qquad\quad{} + \sum_j \bigl(A_j^V(x + \nu_{\ell}/V) -
A_j^V(x)\bigr) V
D_{kj}(t,x + \nu_{\ell}/V)\\
&&\qquad\quad{} + \sum_j \bigl(A_j^V(x + \nu_k/V) -
A_j^V(x)\bigr)VD_{j\ell}(t,x + \nu_k/V) + g_{k\ell}(t,x),
\end{eqnarray*}
where
\begin{eqnarray*}
g_{k \ell}(t,x) &\eqdef& \sum_j V^2 [A_j^V(x + \nu_{\ell}/V +
\nu_k/V)\\
&&\hphantom{\sum_j V^2 [}
{} - A_j^V(x + \nu_{\ell}/V) - A_j^V(x + \nu_k/V) + A_j^V(x)
] \\
&&\hphantom{\sum_j }
{} \times D_j(t,x + \nu_{\ell}/V + \nu_k/V).
\end{eqnarray*}
By (i)\vspace*{1pt} the fact that the second derivative of $A_j^V$ is uniformly
(in $j$ and $x$) bounded and (ii) the bound (\ref{eq:reg1}), the
absolute value of the last term is uniformly (in $t \le T$, $x, k$
and $\ell$) bounded by some $C_3>0$.

As we did for both $v(t,x)$ and $D_k(t,x)$, we change perspective by
viewing the above as a linear system with state space $\{k,\ell,x\}
\in[1,\ldots,M] \times[1,\ldots,M]\times\LL= \Gamma_2$, where we
again put an ordering on $\Gamma_2$ and consider~$\R^{\Gamma_2}$
defined similarly to $\R^{\Gamma_1}$. Also similarly to before, we
note that only a~finite number of the $D_{k,\ell}(t,x)$ are changing
in time. For $b = \{k,\ell,x\} \in\Gamma_2$, we see that $D_b(t)$
satisfies
%
%
\begin{equation}\label{eq:final?}
D_b'(t) = -S_b V D_b(t) + VR_b \cdot D(t) + r_b\cdot D(t) + g_b(t),
\end{equation}
where $D_b(t), D(t), S_b$, $R_b$ and $r_b$ are defined similarly as
before and where we retain the necessary inequalities: for $v \in
\R^{\Gamma_2}$,
%
%
\begin{eqnarray} \label{eq:newbounds}
|R_b \cdot v| &=& \biggl|\sum_j A^V_{j}(b_3) v_{\{b_1,b_2,b_3 +
\nu_j/V\}} \biggr|
\le S_b \|v\|_{\infty}, \nonumber\\[-8pt]\\[-8pt]
|r_b \cdot v| &\le& 2KM \| v\|_{\infty}.
\nonumber
\end{eqnarray}
The rest of the proof is similar to the proof that the $D_k(t,x)$
are uniformly bounded. There is a $b^1 \in\Gamma_2$ and a $t_1 \in
(0,T]$ for which $|D_{b^1}(t)| = \|D(t)\|_{\infty}$ for all $t \in
[0,t_1]$. Taking the derivative of $D_{b^1}(t)^2$ while using
(\ref{eq:final?}), integrating,\vadjust{\goodbreak} and using the bounds
(\ref{eq:newbounds}), we have that for this $b^1$ and any $t \in
[0,t_1]$,
\begin{eqnarray*}
D_{b^1}(t)^2 &=& D_{b^1}(0)^2 + 2\int_0^t g_{b^1}(s)D_{b^1}(s) \,ds -
2\int_0^{t} S_{b^1}VD_{b^1}(s)^2\,ds \\
&&{} + 2\int_0^{t}
VR_{b^{1}}\cdot D(s) D_{b^{1}}(s)\,ds + 2\int_0^{t} r_{b^{1}} \cdot
D(s) D_{b^{1}}(s)\,ds \\
& \le& D_b(0)^2 + 2C_3t + (4KM + 2C_3)\int_0^{t}
\|D(s)\|_{\infty}^2 \,ds,
\end{eqnarray*}
where we used the inequality $x \le1 + x^2$ on the term
$D_{b^{1}}(s)$ in the first integral above. Therefore, for $t \le
t_1$
\[
\|D(t)\|^2_{\infty} \le\|D(0)\|^2_{\infty} + 2C_3t + (4KM +
2C_3)\int_0^{t} \|D(s)\|_{\infty}^2 \,ds.
\]
We continue now by choosing a $b_2 \in\Gamma_2$ such that
$|D_{b_2}(t)| = \|D(t)\|_{\infty}$ for all $t \in[t_1,t_2)$, with
$t_1< t_2 \le T$. By similar arguments as above, we have that for $t
\in[t_1,t_2]$,
\begin{eqnarray*}
\|D(t)\|^2 &\le& \|D(t_1)\|^2_{\infty} + 2C_3(t - t_1) +
(4KM + 2C_3) \int_{t_1}^{t} \|D(s)\|_{\infty}^2 \,ds\\
&\le& \|h(0)\|^2_{\infty} + 2C_3t + (4KM + 2C_3)\int_0^{t}
\|h(s)\|_{\infty}^2 \,ds.
\end{eqnarray*}
Continuing in this manner shows that the above inequality holds for
all $t \in[0,T]$ and so a Gronwall inequality gives us that for all
$t \le T$,
\[
\|D(t)\|_{\infty}^2 \le\biggl( \|D(0)\|^2_{\infty} +
\frac{2C_3}{4KM + 2C_3} \biggr) e^{(4KM + 2C_3)T},
\]
which, after taking square roots, is equivalent to (\ref{eq:reg2}).
\end{pf*}
\begin{pf*}{Proof of Lemma \ref{lemma:lemma_t}}
By (\ref{eq:firstDifference}), we have that for any $k \in
[1,\ldots,M]$ and $x \in\LL$,
\begin{eqnarray*}
D_k(t,x) &=& D_k(0,x) + \sum_j A_j^V(x)\int_0^t D_{kj}(s,x)\,ds \\
&&{} +\sum_j \bigl(A_j^V(x + \nu_k/V) - A_j^V(x)\bigr)V \int_0^t
D_j(s,x+\nu_k/v)\,ds.
\end{eqnarray*}
The proof is now immediate in light of Lemma \ref{lemma:reg}.\vadjust{\goodbreak}
\end{pf*}

\section{Examples}
\label{sec:examples}

\begin{example}
Consider the case of an irreversible isomerization of one molecule
into another. We denote by $A$ the molecule undergoing the
isomerization and $B$ the target molecule. We assume that the rate
constant associated with this reaction is $1$. The pictorial
representation for this system is simply
\[
A \stackrel{1}{\to} B.
\]
Letting $X(t)$ denote the number of $A$ molecules at time $t\ge0$,
$X(t)$ satisfies
\[
X(t) = X(0) - Y\biggl(\int_0^t X(s) \,ds\biggr).
\]
Supposing that we start with $V = 10\mbox{,}000$ molecules, we approximate
the distribution of $X(1)$ using $200\mbox{,}000$ sample paths constructed
using the Gillespie algorithm, which produces statistically exact
sample paths, Euler tau-leaping with a step-size of $1/20$
and\vspace*{1pt}
midpoint tau-leaping with a step-size of $1/20$. Note that in this
case $1/20 = 1/V^{0.325}$, and so $\beta= 0.325$. The computational
results are presented in Figure \ref{fig:iso}, which demonstrate
the stronger convergence rate of midpoint tau-leaping as compared to
Euler tau-leaping.

%
%
\begin{figure}[t]

\includegraphics{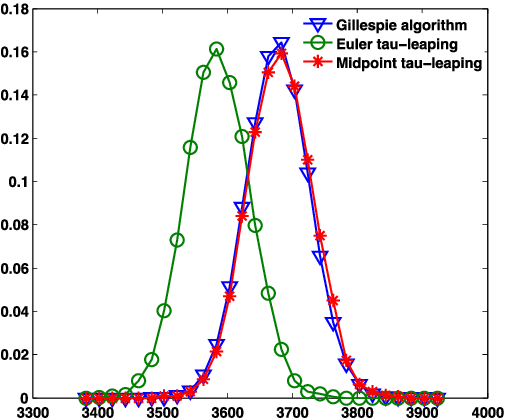}

\caption{Relative frequency of $X(1)$ from $200\mbox{,}000$ sample paths
constructed using \textup{(i)} Gillespie's algorithm, blue line
$\nabla$ marker, \textup{(ii)} Euler tau-leaping, green line, $\bigcirc$ marker,
and \textup{(iii)}~midpoint tau-leaping, red line, $*$ marker. The
approximated distribution generated via midpoint tau-leaping is
clearly closer to the exact distribution than that of Euler
tau-leaping.}
\label{fig:iso}
\end{figure}

It is simple to show that $X(1)$ is a binomial$(n,p)$ random
variable with parameters $n=10\mbox{,}000$ and $p=1/e$. Therefore, $\E
X(1) = 10\mbox{,}000/e \approx3678.8$. The estimated means produced from
the 200,000 sample paths of Euler tau-leaping and midpoint
tau-leaping were $3585.4$ and $3681.4$, respectively. Solving for
$\Err(t)$ of (\ref{eq:det}) for this example yields $\Err(t) =
(1/2)e^{-t}t$. Theorem \ref{thm:weak_tau_exact} therefore estimates
that Euler tau-leaping should produce a mean $(1/2)e^{-1} 10\mbox
{,}000^{1 -
0.325} \approx92.2$ smaller than the actual mean, which is in
agreement with $3678.8 - 3585.4 = 93.4$. Solving for $\Err_1(t)$ of
(\ref{err1}) for this example yields $\Err_1(t) = (1/6)te^{-t}$.
Theorem \ref{thm:weak_mdpt_exact} therefore estimates that midpoint
tau-leaping should produce a mean $(1/6)e^{-1}10\mbox{,}000^{1 -
2*0.325} =
4.62$ smaller than the actual mean, which is in agreement with
$3678.8 - 3681.4 = -2.6$.
\end{example}
\begin{example}
We now consider a simple Lotka--Volterra predator--prey model.
Letting $A$ and $B$ represent the prey and predators, respectively,
in a~given environment we suppose (i) prey reproduce at a certain
rate, (ii)~interactions between predators and prey benefit the
predator while hurting the prey, and (iii) predators die at a
certain rate. One possible model for this system is
\[
A \stackrel{2}{\rightarrow} 2A,\qquad
A+B\stackrel{0.002}{\rightarrow} 2B,\qquad
B\stackrel{2}{\rightarrow}\varnothing,
\]
where a choice of rate constants has been made. Letting $X(t) \in
\Z^2_{\ge0}$ be such that $X_1(t)$ and $X_2(t)$ represent the
numbers of prey and predators at time $t>0$, respectively, $X(t)$
satisfies
%
%
\begin{eqnarray}\label{eq:LV}
X(t) &=& X(0) + Y_1\biggl(\int_0^t
2X_1(s)\,ds\biggr)\left[\matrix{
1\cr
0}
\right]\nonumber\\
&&{} + Y_2\biggl(\int_0^t
0.002X_1(s)X_2(s) \,ds\biggr)
\left[\matrix{
-1\cr
1}
\right]\\
&&{} + Y_3\biggl(\int_0^t
2X_2(s) \,ds\biggr)
\left[\matrix{0\cr
-1}\right].\nonumber
\end{eqnarray}
We take $X(0) = [1\mbox{,}000,1\mbox{,}000]^T$, and so $V = 1\mbox
{,}000$ for our model.
Lotka--Volterra models are famous for producing periodic solutions;
this behavior is demonstrated in Figure \ref{fig:LV1}.

%
%
\begin{figure}

\includegraphics{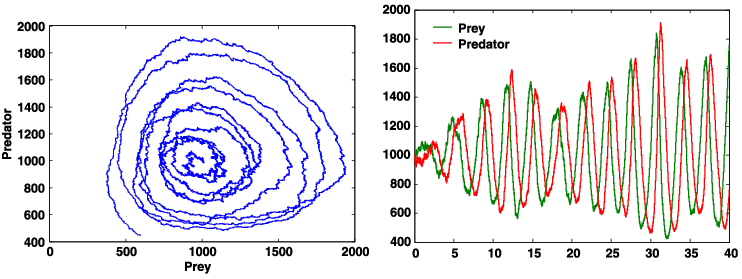}

\caption{Oscillations in a predator--prey model. In the left image
we see the numbers of predators versus the number of prey for a
single realization of the system (\protect\ref{eq:LV}). In the right
image we see the time-series of the numbers of predators and
prey for a single realization of (\protect\ref{eq:LV}).}
\label{fig:LV1}
\end{figure}

%
%
\begin{figure}[t]

\includegraphics{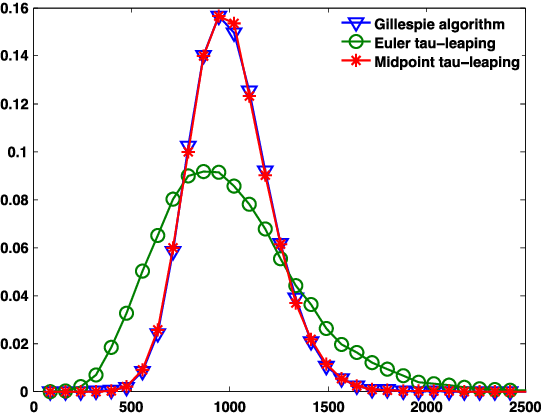}

\caption{Relative frequency of $X_2(10)$ from $30\mbox{,}000$ sample
paths constructed using \textup{(i)}~Gillespie's algorithm, blue line $\nabla$
marker, \textup{(ii)} Euler tau-leaping, green line, $\bigcirc$ marker, and
\textup{(iii)} midpoint tau-leaping, red line, $*$ marker. The
approximated distribution generated via midpoint tau-leaping is
clearly closer to the exact distribution than that of Euler
tau-leaping.}
\label{fig:LVdist}
\vspace*{-2pt}
\end{figure}

We approximate the distribution of $X_2(10)$ using $30\mbox{,}000$ sample
paths constructed using the Gillespie algorithm, Euler tau-leaping
with a step-size of $1/20$ and midpoint tau-leaping with a
step-size of $1/20$. Note\vspace*{1pt} that in this case $1/20 = 1/V^{0.434}$,
and so $\beta= 0.434$. The computational results are presented in
Figure \ref{fig:LVdist}, which again demonstrate the stronger
convergence rate of midpoint tau-leaping as compared to Euler
tau-leaping.
\end{example}

\begin{appendix}\label{sec:app}

\section*{Appendix}

\setcounter{theorem}{0}
\begin{lemma}\label{mgincrem}
Let $M$ be a $\{\mathcal{F}_t\}$-martingale, $R$ be bounded and
$\{\mathcal{F}_t\}$-adapted, and let $h>0$. Then for $\eta(t)\equiv[t/h]h$,
\[
\hat{M}(t)= \int_0^tR \circ\eta(s) \bigl(M(s) - M \circ\eta(s) \bigr) \,ds
+ R \circ\eta(t)\bigl(M(t)-M\circ\eta(t)\bigr)\bigl(\eta(t)+h-t\bigr)
\]
is an $\{\mathcal{F}_t\}$-martingale and
%
%
\setcounter{equation}{0}
\begin{equation}\label{mhatvar}
[\hat{M}]_t=\int_0^t\bigl(R\circ\eta(r)\bigr)^2\bigl(\eta(r)+
h-r\bigr)^2\,d[M]_r.
\end{equation}
If $M$ is $\mathbb{R}^d$-valued and $R$ is $ \mathbb{M}^{m\times
d}$-valued, then the quadratic covariation matrix is
\[
[\hat{M}]_t=\int_0^t\bigl(\eta(r)+h-r\bigr)^2R\circ\eta
(r)\,d[M]_rR^T\circ\eta(r).
\]
\end{lemma}
\begin{pf}
For $t<T-h$,
\begin{eqnarray*}
\E[\hat{M}(T)|\mathcal{F}_t] &=& \E\biggl[\int_0^TR\circ\eta
(s)\bigl(M(s)-M\circ\eta(s)\bigr)\,ds\big|\mathcal{F}_
t\biggr]\\[-1pt]
&=& \E\biggl[\int_0^{\eta(t)+h}R\circ\eta(s)\bigl(M(s)-M\circ\eta
(s)\bigr)\,ds\big|\mathcal{F}_t\biggr]\\[-1pt]
&=& \int_0^tR\circ\eta(s)\bigl(M(s)-M\circ\eta(s)\bigr)\,ds\\[-1pt]
&&{}+R\circ\eta(t)\bigl(M
(t)-M\circ\eta(t)\bigr)\bigl(\eta(t)+h-t\bigr).
\end{eqnarray*}
The case of $T-h\le t < T$ is similar. $[\hat{M}]$ is just the
quadratic variation of the second term on the right, and noting that
$\hat{M}$ is continuous at $t=kh$ for all $k=0,1,2\ldots,$
(\ref{mhatvar}) follows.
\end{pf}

For completeness, we include a statement of the martingale central
limit theorem (see \cite{Kurtz86} for more details).\vadjust{\goodbreak}
\begin{lemma}\label{lem:MCLT}
Let $\{M_n\}$ be a sequence of $\R^d$-valued martingales with $M_n(0)
= 0$. Suppose
\[
\lim_{n \to\infty}\E\Bigl[ \sup_{s \le t}|M_n(s) - M_n(s-)|\Bigr] = 0
\]
and
\[
[M_n^i,M_n^j]_t \to c_{i,j}(t)
\]
for all $t\,{>}\,0$ where $C\,{=}\,((c_{i,j}))$ is deterministic and
continuous. Then \mbox{$M_n\,{\Rightarrow}\,M$}, where $M$ is Gaussian with
independent increments and \mbox{$\E[M(t)M(t)^T] = C(t)$}.
\end{lemma}
\end{appendix}

%
%

%
\printaddresses

\end{document}